\documentclass[moor]{informs1}              

\usepackage[ruled,lined,boxed,norelsize]{algorithm2e}
\SetKwBlock{Repeat}{repeat}{end}

\usepackage{tikz}
\usetikzlibrary{calc,patterns,angles,quotes}
\usepackage{soul}
\usepackage{pifont}

\usepackage{amsmath,amsfonts,amssymb}

\usepackage{graphics}
\usepackage{epsfig}

\usepackage{inputenc}
\usepackage[T1]{fontenc}

\usepackage{color}
\usepackage{float}
\usepackage{bm}
\usepackage{enumerate}
\usepackage[noadjust]{cite}
\usepackage{graphicx}
\usepackage{comment}
\usepackage{marginnote}
\usepackage{xcolor}
\usepackage{mathtools}
\DeclarePairedDelimiter\ceil{\lceil}{\rceil}

\def\trt{^{\scriptscriptstyle T}}

\usepackage{natbib}
 \NatBibNumeric
 \def\bibfont{\small}%
 \bibpunct[, ]{[}{]}{,}{n}{}{,}%

\usepackage[colorlinks=true,breaklinks=true,bookmarks=true,urlcolor=blue,
     citecolor=blue,linkcolor=blue,bookmarksopen=false,draft=false]{hyperref}

\def\EMAIL#1{\href{mailto:#1}{#1}}
\def\URL#1{\href{#1}{#1}}         
\renewcommand{\int}{\mathop{\rm int}}

\TheoremsNumberedThrough     

\EquationsNumberedThrough    

\begin{document}
\TITLE{Ghost penalties  in nonconvex constrained optimization: Diminishing stepsizes and iteration complexity}
\ARTICLEAUTHORS{%
\AUTHOR{Francisco Facchinei}
\AFF{Sapienza University of Rome,
Department of Computer, Control, and Management Engineering Antonio Ruberti, Via Ariosto 25, 00185 Rome, \EMAIL{francisco.facchinei@uniroma1.it}, \URL{}}
\AUTHOR{Vyacheslav Kungurtsev}
\AFF{Czech Technical University in Prague, Department of Computer Science, Faculty of Electrical Engineering, \EMAIL{vyacheslav.kungurtsev@fel.cvut.cz}, \URL{}}
\AUTHOR{Lorenzo Lampariello}
\AFF{Roma Tre University,
Department of Business Studies, Via D'Amico 77, 00145 Rome, \EMAIL{lorenzo.lampariello@uniroma3.it}, \URL{}}
\AUTHOR{Gesualdo Scutari}
\AFF{Purdue University, School of Industrial Engineering,  West-Lafayette, IN, USA, \EMAIL{gscutari@purdue.edu}, \URL{}}
} 
\ABSTRACT{We consider nonconvex constrained optimization problems and propose a new approach to the convergence analysis based on penalty functions. We make use of classical penalty functions in an unconventional way, in that penalty functions only enter in the theoretical analysis of convergence while the algorithm itself is penalty-free. Based on this idea, we are able to establish several new results, including the first general analysis for diminishing stepsize methods in nonconvex, constrained optimization, showing convergence to generalized stationary points, and a complexity study for SQP-type algorithms.

}%
\MSCCLASS{90C30, 90C60, 65K05}
\ORMSCLASS{Primary: programming: nonlinear  algorithms; secondary: analysis of algorithms: computational complexity }
\HISTORY{}
\maketitle
%
\section{Introduction}\label{Sec:Intro}
We consider the nonconvex constrained optimization problem
\begin{equation}\label{eq:startpro}
\begin{array}{cl}
\underset{x}{\mbox{minimize}} & f(x) \\
\mbox{s.t.} & g(x) \le 0\\[5pt]
& x \in K,
\end{array}\tag{P}
\end{equation}
where $K \subseteq \mathbb R^n$ is a nonempty closed and convex set, and $f : \mathbb R^n \to \mathbb R $
and $g : \mathbb R^n \to \mathbb R^{m} $ are C$^{1,1}$ (i.e., continuously differentiable with locally Lipschitz
gradients) functions on an open set containing $K$.

Penalty functions (differentiable or nondifferentiable, exact or sequential) are part of the folklore in optimization
and have been widely used in analyzing optimality conditions, stability and sensitivity
properties and in developing solution methods. In this paper we put forward a new use of penalty functions in
the design of algorithms for the solution of \eqref{eq:startpro}. In particular, we consider the classical
nondifferentiable penalty function
\[ W(x; \varepsilon) \triangleq f(x) + \frac{1}{\varepsilon} \max_i\{g_i(x)_+\},
 \]
 where $\alpha_+ \triangleq  \max \{0,\alpha\}$ and $\varepsilon$ is a positive penalty parameter.
 We propose a novel use of  $W(x; \varepsilon)$ in that, contrary to usual penalty algorithms, the penalty
 function {\em only enters in the theoretical analysis} of convergence while {\em the algorithm itself is penalty-free,} whence the term {\em ghost} penalty.
 We establish (subsequential)
 convergence to ``generalized stationary points" under essentially no assumptions beyond the C$^{1,1}$
 requirement on the problem functions. In particular, we assume neither feasibility of the problem nor any
 constraint qualification and therefore (subsequential) convergence to generalized stationary point is the natural target for a well-behaved algorithm, see for example \cite{birgin2016evaluation,burke1989sequential,burke1992robust,  burke1989robust,burke2013epi,cartis2011evaluation, cartis2013evaluation,cartis2017corrigendum,cartis2019evaluation, cartis2018optimality,el1999global,Facch97,liu2000robust, liu2004robust,martinez2017high, yuan1995convergence}.
Once our main convergence result has been established, the role played by further classical  assumptions, like feasibility or constraint qualifications, for example,  is easily understood and can be taken into account in a straightforward way.

Our framework is of a generalized Sequential Quadratic Programming (SQP)-type. At each iteration $x^\nu$ we generate a search direction $d(x^\nu)$ by solving a strongly convex optimization subproblem constructed along the lines discussed in the seminal papers \cite{burke1989sequential,burke1989robust} and also taking into account the developments in \cite{facchinei2017feasible,scutari2014decomposition}. The direction-finding subproblem reads as follows:
\begin{equation}\label{eq:approx problem intro}
   \underset{d}{\mbox{minimize}} \; \tilde f (d;x^\nu) \quad \text{s.t.} \quad
   \tilde g (d;x^\nu) \le \kappa(x^\nu) \, e, \; \|d\|_\infty \le \beta, \; d \in K - x^\nu,
\end{equation}
where $\tilde f (\bullet ;x^\nu)$ and $\tilde g_i (\bullet; x^\nu)$ are strongly convex and convex
approximations of $f$ and $g$, respectively, $\kappa(x^\nu) \in \mathbb R$ is nonnegative and defined to
make the subproblem feasible, $\beta$ is a user-chosen positive constant, and $e \in \mathbb R^m$ is the
vector with all components being one. The classical SQP subproblem is a particular instance of \eqref{eq:approx
problem intro}, when $\tilde g$ is just a linearization of $g$ and $\tilde f$ a positive definite quadratic
approximation of $f$. Our more general approach leaves room for much flexibility in tailoring the direction
finding subproblem to the problem at hand and to exploit any available specific structure in \eqref{eq:startpro},
see Section \ref{Sec:dirfind} for more details.

A step $\gamma ^\nu$ is then taken along this direction so that
\begin{equation}\label{eq:method}
x^{\nu+1} \,=\, x^\nu + \gamma^\nu d(x^\nu).
\end{equation}
 We consider both classical Diminishing Stepsize Methods (DSMs) wherein $\gamma^\nu$ is a positive stepsize
 such that
\begin{equation}\label{eq:gamma}
\lim_{\nu\to \infty} \gamma^\nu = 0\qquad \mbox{\rm and}\qquad \sum_{\nu=0}^\infty \gamma^\nu\, =\,
\infty,
\end{equation}
and stepsize selection rules where $\gamma^\nu$ is chosen in a more problem-tailored way, typically fixed for
 at least a   subsequence of iterations if not for the entire sequence.
Note that, while  the algorithms in this paper generate the search direction in a (generalized) SQP-fashion, they
differ markedly from classical SQP methods in the way they select stepsizes. Indeed, SQP methods traditionally
have adopted effective globalization procedures based e.g.  on line-search or trust-region strategies. Here,
instead,  we mostly study different techniques that may be useful, for example, in very large-scale or
distributed settings and that, in addition, allow us to perform an iteration complexity analysis. Given the
effectiveness of the SQP approach in handling nonconvex constraints, our results hopefully provide an
alternative to expand the applicability of SQP-like methods.

It may be interesting to mention at this juncture the two papers \cite{auslender2010moving}  and
\cite{bolte2016majorization} where, in the context of an extended SQP-like scheme, a  stepsize of one is always taken,
thus also foregoing line-search, trust-region or other standard globalization techniques.
The possibility to use a fixed (large) stepsize derives
from the fact that the methods in \cite{auslender2010moving}  and \cite{bolte2016majorization}  are {\em
feasible} methods and the surrogate for the objective function is always an {\em
upper} convex approximation; in fact, the algorithms analyzed in these two interesting works belong to the class of SQP-type majorization-minimization schemes. The essential ingredient in developing such methods is the ability to build approximations that majorize both the constraints and the objective function; we will discuss the consequences of this setting in the context of our scheme in Sections \ref{Sec:complexity} and \ref{sec:bounded}.

Some other related papers in the contemporary literature include \cite{auslender2013extended,cartis2011evaluation,fletcher19851,solodov2004sequential}, where penalty algorithms are analyzed.
The methods discussed in these works aim at minimizing directly the penalty function by resorting to composite optimization approaches; this results in a double-loop structure whereby the penalty function is (approximatively) minimized for a certain value of the penalty parameter, and then the penalty parameter is possibly updated.
Although there are some similarities in the analysis, the algorithms presented in our paper rely on a pure ``SQP-like'' subproblem and
penalty parameters enter only in the theoretical analysis, and not in the subproblem definition.

Based on our ghost penalty approach, our main contributions are

 (a) the first (subsequential) convergence result for a wide class of DSMs for general  nonconvex constrained optimization problems;

 (b) {\em iteration complexity results} for some suitable choices of the stepsize $\gamma^\nu$ leading, among other things, to the first iteration complexity analysis for SQP-type methods in a general setting.
\smallskip

 \noindent
 The results related to (a) considerably widen the scope of applicability of DSMs. DSMs are part of the core techniques in optimization and their advantages and disadvantages are well known, see for example \cite{bertsekas2015convex,polyak1987introduction,shor2012minimization}.
However, DSMs are not yet fully understood when it comes to  nonconvex problems. Indeed, the very recent paper \cite{davis2019stochastic} is the first study to establish convergence results for DSMs applied to unconstrained, nonsmooth, and nonconvex problems. The results in the present paper, therefore, close a surprising gap in the literature, since DSMs have never been proved to lead to convergence when addressing problems with nonconvex constraints except in some specialized settings where feasibility of the iterates can be maintained throughout the optimization process and constraint qualifications are assumed \cite{facchinei2017feasible,scutari2017parallel}.
We show that every limit point of the sequence $\{x^\nu\}$ produced according to \eqref{eq:method} and \eqref{eq:gamma} is a generalized stationary solution for \eqref{eq:startpro}. By generalized stationary, we intend a point that can be:  an infeasible stationary solution of the violation-of-the-constraint problem associated to \eqref{eq:startpro}, a Fritz-John or a KKT point of \eqref{eq:startpro}. As mentioned above, this is the natural target of an algorithm for constrained optimization  when neither  blanket assumptions about feasibility of \eqref{eq:startpro} are made, nor constraint qualifications are assumed to hold.
Many consider DSMs a ``necessary evil" and nevertheless they are
currently used in many settings, for example in parallel and distributed optimization, in stochastic optimization,  in multi-agent settings, in incremental methods, and whenever the computation of the objective function is very expensive or the problem is affected by noise, see for example
  \cite{bertsekas1989parallel,bertsekas2011incremental,bertsekas2015convex,
  bottou2018optimization,davis2019stochasticA,davis2019stochastic,
 daneshmand2015hybrid,daneshmand2018decentralized,dean2012large,di2016next,
 facchinei2015parallel,facchinei2017feasible,Nedich_et_al_2010,gupta2014training,kingma2014adam,
 nemirovski2009robust,ng2007joint,
  polyak1987introduction,scutari2014decomposition,scutari2017parallel,scutari2017parallelbis,
 tatarenko2017non,wang2017stochastic,zeng2018nonconvex}. In some cases it might be hoped that soon more effective alternatives will be found, in other cases alternatives are harder to anticipate.
In any case, as DSMs evolve to deal with new classes of problems of contemporary interest, the need to tackle nonconvex constraints and to relax the conditions needed to analyze convergence emerges.
The developments in this paper, dealing with a standard nonconvex optimization problem,  are a first step in this direction and will hopefully pave the way for further advancements in the more challenging settings mentioned above.

Results indicated in (b) add to a thus-far sparse, but thriving literature that just recently began appearing on the topic of
complexity  analysis for nonconvex optimization problems.
Disregarding classical results on the gradient method, see e.g. \cite{nesterov2013introductory},
this chapter was opened by the largely ignored  Vavasis' paper \cite{vavasis1993black}, but gained momentum only with  Nesterov and Polyak's work \cite{nesterov2006cubic} on a cubic
regularization method for the unconstrained minimization of a nonconvex, smooth function.
An excellent review of results in this field is contained in \cite{cartis2018optimality}, to which we refer the interested reader for a broader view on the subject. Here we only discuss results on algorithms for nonconvex, inequality constrained problems aimed at locating {\em generalized} stationary points  using first-order information, similarly to what we do in the present paper.

By using our ghost penalty approach we are able to establish that ${\mathcal{O}}(\delta^{-4})$ iterations are needed at worst to find a $\delta-$approximate generalized stationary point, \textcolor{black}{see Definition \ref{def:delta approx};  this definition of $\delta-$approximate generalized stationarity relaxes the notion of an exact generalized stationary point and parallels similar developments in \cite{birgin2016evaluation,cartis2017corrigendum, cartis2018optimality, cartis2019evaluation}. In line with what discussed above, we remark that
 our notion of (approximate) stationarity is naturally broader than the more standard (approximate) KKT conditions. }
The bound of ${\mathcal{O}}(\delta^{-4})$ can be reduced to ${\mathcal{O}}(\delta^{-3})$ if a feasible point is known in advance and to ${\mathcal{O}}(\delta^{-2})$ if some further conditions are met (see Section \ref{Sec:complexity}  for details).
As far as we are aware of, these are the  first iteration complexity results for SQP-type methods
 \textcolor{black}{{\em in a general setting}, i.e. without assuming, for example, feasibility of iterates, see the discussion about \cite{auslender2010moving} and \cite{bolte2016majorization}
below.} Indeed, with the exception \cite{cartis2011evaluation},  all other results for general nonconvex, constrained problems obtained so far in the literature are based on Phase I - Phase II  type methods wherein an almost feasible point is first sought and then a second phase is started.
More specifically,
in \cite{cartis2011evaluation} a penalty approach is shown to take either ${\mathcal{O}}(\delta^{-2})$ or ${\mathcal{O}}(\delta^{-5})$ iterations to reach an approximate generalized stationary point, depending on the behavior of the penalty parameter during the minimization process. Cartis, Gould and Toint also describe in \cite{cartis2013evaluation} a Phase I - Phase II cubic regularization method, possibly using Hessian information, for the solution of equality constrained problems and show that the number of iterations needed to reach an approximate generalized stationary point varies between ${\mathcal{O}}(\delta^{-3/2})$ and ${\mathcal{O}}(\delta^{-2})$ depending on certain algorithmic choices.
Building on the algorithm in \cite{cartis2013evaluation}, and using a different definition of approximate generalized stationary point,  Birgin et al. show in \cite{birgin2016evaluation} that a Phase I - Phase II algorithm takes between ${\mathcal{O}}(\delta^{-3})$ and ${\mathcal{O}}(\delta^{-5})$ iterations to declare a point approximate generalized stationary, according to the choice of an algorithmic parameter. Finally, Cartis, Gould, and Toint establish in  \cite{cartis2017corrigendum} a bound of ${\mathcal{O}}(\delta^{-2})$, once again for a Phase I - Phase II method and using first-order information only.
A  detailed comparison of all these results is not straightforward, because of the many subtleties involved, and we defer a more detailed discussion  on this issue to Remark \ref{rem:compl} in Section \ref{Sec:complexity}.
We conclude mentioning, once again, the SQP-like majorization-minimization methods proposed in  \cite{auslender2010moving} and \cite{bolte2016majorization}. \textcolor{black}{Differently to what we propose here, these two papers discuss only {\em feasible-type} methods and assume standard constraint qualifications.  In this framework,  interesting results are obtained concerning the \emph{convergence rate} to zero of the {\em distance of the point generated by the method to a KKT solution} (as opposed to the more algorithmically oriented results in the papers discussed above, where bounds are obtained on the number of iterations needed to satisfy a given stopping criterion).}
Note that the distinction of iteration complexity and convergence rate is a subtle and sometimes blurred
one. In a nutshell, the difference is that when we talk about complexity we are assuming that the constants appearing in the complexity bound are  conceptually known a priori (e.g. Lipschitz constants), while in the case of a convergence rate the bounds include constants that are possibly unknown in advance (e.g. the diameter of the region that contains all iterates).
In \cite{auslender2010moving}, linear convergence of the sequence of iterates to the optimal solution is obtained for  strongly convex problems.
 The  more general results in \cite{bolte2016majorization} dispense with convexity by assuming the
 Kurdyka-\L{}ojasiewicz property and obtaining a convergence rate that depends on the
 \L{}ojasiewicz exponent.



The paper is organized as follows. In Section \ref{sec:stationary} we introduce some mathematical preliminaries and, in particular, the appropriate definition of a generalized stationary point for nonconvex, constrained problems. In Section \ref{Sec:dirfind} we discuss in detail the direction finding subproblem and introduce some assumptions that will be used to establish convergence. In Section \ref{Sec:convergence} we show convergence to generalized stationary points for DSMs, while in Section \ref{Sec:complexity} we perform the iteration complexity analysis. We finish in Section~\ref{sec:bounded} with a discussion on the boundedness of the sequence of iterates.

\section{Generalized Stationary Points}\label{sec:stationary}
We consider  Problem \eqref{eq:startpro}, under the blanket assumptions indicated in the Introduction, and denote the feasible set of \eqref{eq:startpro} by
\[
{\mathcal{X}} \triangleq \left\{d \in \mathbb R^n \, : \, g(x) \le 0, \, d \in K\right\}.
\]
Note that we do not assume that problem \eqref{eq:startpro} is feasible, let alone that it has a solution. Therefore, we aim at deriving convergence results for both feasible and infeasible problems, in some suitable sense.

A general constrained problem \eqref{eq:startpro} can be viewed as a combination of two problems: (i) the feasibility one, i.e., the problem of finding a feasible point; and (ii) the problem of finding a local minimum point of the objective function over the feasible set. Even just the former problem is a hard one, since it essentially requires computing a global minimum of the generally nonconvex function expressing the violation of the constraints. Consistently, we design our algorithm to converge to stationary solutions in a generalized sense, that is to points that either are stationary for \eqref{eq:startpro} or are infeasible and stationary for the following violation-of-the-constraints optimization problem:
\begin{equation}\label{eq:feaspro}
\begin{array}{cl}
\underset{x}{\mbox{minimize}} & \displaystyle \max_i \{g_i(x)_+\},  \\
& x \in K
\end{array}
\end{equation}
where, we recall,  $\alpha_+ \triangleq \max\left\{0, \alpha\right\}$ for all $\alpha \in \mathbb R$.
Let
\begin{equation*}\label{eq:norm}
M_1(x) \triangleq \left\{\xi \, \left| \right. \, \xi \in N_{\mathbb R^m_-}(g(x)), \, 0 \in \nabla f(x) + \nabla g(x) \xi + N_K(x) \right\}
\end{equation*}
and
\begin{equation*}\label{eq:abn}
M_0(x) \triangleq \left\{\xi \, \left| \right. \, \xi \in N_{\mathbb R^m_-}(g(x) - \max_i \{g_i(x)_+\} e), \, 0 \in \nabla g(x) \xi + N_K(x) \right\},
\end{equation*}
where $N_{\mathbb R^m_-}(y)$ and $N_K(x)$ are the classical normal cones to the convex sets $\mathbb R^m_-$ and $K$ at $y$ and $x$, respectively, $\nabla f$ is the gradient of $f$ and $\nabla g$ is the transposed Jacobian of $g$. If  $g(x) \le 0$, the condition $\xi \in N_{\mathbb R^m_-}(g(x))$ can be more familiarly rewritten as $\xi_i \geq 0,\, \xi_i g_i(x) =0$ for all $i$ (a similar reasoning applies to the normal cone expression in the definition of $M_0(x)$).  We note explicitly that if $x$ is not feasible, the set $M_1(x)$ is empty.
Let  $\hat x$  be a local minimum point of \eqref{eq:startpro}, then it is well-known that either $M_1(\hat x) \neq \emptyset,$ (the point is a KKT point) or $M_0(\hat x) \neq \{0\}$ (the point is a Fritz-John point), or both. On the contrary, it is classical to show that if $\hat x \in K$ is not feasible, i.e. if  $g_i(\hat x) > 0$ for at least one index $i \in \{1, \ldots, m\}$, in view of the regularity of the functions involved, then the stationarity condition for problem \eqref{eq:feaspro},
\begin{equation}\label{eq:statfeaspro}
0 \in \partial \max_i\{g_i(\hat x)_+\} + N_K(\hat x),
\end{equation}
is equivalent to $M_0(\hat x) \neq \{0\}$. Hence, the (generalized) stationarity criteria for the original problem \eqref{eq:startpro} can naturally be specified by using the  sets $M_1$ and $M_0,$ as detailed in Definition \ref{df:stationarity}.
\begin{definition}\label{df:stationarity}
A point $\hat x \in K$ is, for problem \eqref{eq:startpro},
\begin{enumerate}
\item[$\bullet$] a KKT solution if $g(\hat x) \le 0$ and $M_1(\hat x) \neq \emptyset;$
\item[$\bullet$] a Fritz-John (FJ) solution if $g(\hat x) \le 0$ and $M_0(\hat x) \neq \{0\};$
\item[$\bullet$] an External Stationary (ES) solution if $g_i(\hat x) > 0$ for some  $i \in \{1, \ldots, m\}$ and $M_0(\hat x) \neq \{0\}.$
\end{enumerate}
We call $\hat x \in K$ a generalized stationary solution of \eqref{eq:startpro} if any of these cases occurs.
\end{definition}

Since we did not make any regularity or feasibility assumptions on problem \eqref{eq:startpro}, finding a generalized stationary solution in the sense just described is
the appropriate requirement for a solution algorithm; we show that our method does converge to generalized stationary points as defined above.
It also turns out that, under classical regularity conditions, our algorithm actually converges to KKT points. The constraint qualification (CQ) we use is the Mangasarian-Fromovitz one, suitably extended to (possibly) infeasible points.
\begin{definition}\label{df:emfcq}
We say that the extended Mangasarian-Fromovitz Constraint Qualification (eMFCQ) holds at $\hat x \in K$ if
\begin{equation*}\label{eq:emfcq}
M_0(\hat x) = \{0\}.
\end{equation*}
\end{definition}
If $\hat x$ is feasible and $K = \mathbb R^n$, this condition reduces to the classical MFCQ and in turn, whenever the constraints are convex, it is well-known that the MFCQ is equivalent to Slater's CQ, i.e. to the existence of a point $\tilde x$ such that $g(\tilde x) < 0$.
The eMFCQ is rather standard and its definition goes back at least to
\cite{di1988exactness,di1989exact}, having its roots in \cite{rockafellar1982lagrange}; since its introduction it has been used rather often, especially in the analysis of penalty and SQP algorithms, because it arises quite naturally in these contexts.
By using \cite[Motzkin's theorem of alternative 2.5.2]{craven2012mathematical}, we see that the eMFCQ holds at $\hat x \in K$ if and only if
 \begin{equation}\label{eq:EMFCQ1}
 \exists \, \hat d \,\in\, T_K(\hat x):\; \nabla g_i(\hat x)\trt \hat d < 0, \quad \forall i: g_i(\hat x) = \max_j\{ g_j(\hat x)_+\}.
 \end{equation}
 Since $K$ is convex, simple continuity arguments show that the latter condition is equivalent to
  \begin{equation}\label{eq:EMFCQ2}
 \exists \, \tilde x \,\in\, K:\; \nabla g_i(\hat x)\trt (\tilde x - \hat x) < 0, \quad \forall i: g_i(\hat x) = \max_j\{ g_j(\hat x)_+\}.
 \end{equation}
 We state below a result that extends a standard property  of the MFCQ for feasible points.
\begin{proposition}\label{th:cqlocal}
If the eMFCQ holds at $\hat x \in K$, then there exists a neighborhood $\mathcal V$ of $\hat x$ such that, for every $x \in K \cap \mathcal V$, the eMFCQ is satisfied.
\end{proposition}
\proof{Proof.}
 If $\hat x \in K$ is feasible, this is a classical result. If $\hat x \in K$ is not feasible, the condition $M_0(\hat x) = \{0\}$ implies that
$\hat x$ is not a stationary point for the feasibility problem \eqref{eq:feaspro}, i.e. $0\not \in \partial  \max_i \{g_i(\hat x)_+\} + N_K(\hat x)$. The assertion then easily follows from the outer semicontinuity and local boundedness of the subdifferential mapping $\partial \max_i \{g_i(\bullet)_+\}$ and by (see \cite[Proposition 6.6]{RockWets98}) the outer semicontinuity relative to $K$ of the set valued mapping $N_K$ (see \cite{RockWets98} for the definition of outer semicontinuity). \hfill \Halmos\endproof

\section{Direction Finding Subproblem}\label{Sec:dirfind}
At each iteration of our algorithm we move from the current iteration $x^\nu$ along the direction $d(x^\nu)$ with a stepsize $\gamma^\nu$, see \eqref{eq:method}. While the stepsize is chosen according to several rules to be discussed in the following sections, the direction $d(x^\nu)$ is the solution of the strongly convex subproblem \eqref{eq:approx problem intro},  briefly described in the Introduction, that we repeat here for the reader's convenience.

Given a (base) point $x \in K$ (which will actually  be the current iterate $x^\nu$ in the algorithm), $d(x)$ is the unique solution of the following strongly convex optimization problem:
\begin{equation}\label{eq:p_k}
\begin{array}{cl}
\underset{d}{\mbox{minimize}} &  \tilde f(d;x)\\
\mbox{s.t.} & \tilde g(d; x) \le \kappa(x) e\\[5pt]
& \|d\|_\infty \le \beta,\\[5pt]
& d \in K - x
\end{array}\tag{P$_{x}$}
\end{equation}
where $e \in \mathbb R^m$ is the vector with all components being one and $\beta$ is a user-chosen positive constant. Moreover, $\tilde f$ is a strongly convex surrogate of the original objective function $f$, while $\tilde g$ is a convex surrogate of the original constraints $g$ (see Assumption A below for the conditions these surrogates must obey). Finally, following \cite{burke1989sequential}, the quantity $\kappa(x)$ in the surrogate constraints, which serves to suitably enlarge the feasible set of the subproblem in order to ensure it is always nonempty, is defined, for every $x \in K$, as follows:
\begin{equation}\label{eq:cap}
\kappa(x) \triangleq (1 - \lambda) \max_i \left\{g_i(x)_+\right\} + \lambda \min_d \left\{\max_i  \left\{\tilde g_i(d; x)_+ \right\} \, | \, \|d\|_\infty \le \rho, \, d \in K - x\right\},
\end{equation}
with $\lambda \in (0,1)$ and  $\rho \in (0, \beta)$.
Note that \eqref{eq:cap} requires the computation of the optimal value of the convex problem
\begin{equation}\label{eq:capoptpro}
\min_d \left\{\max_i  \left\{\tilde g_i(d; x)_+ \right\}\, | \, \|d\|_\infty \le \rho, \, d \in K - x\right\}
\end{equation}
that always has an optimal solution because the feasible set is nonempty and compact. Note also that if $x$ is feasible for \eqref{eq:startpro}, $\kappa(x) = 0$.  The additional constraint $\|d\|_\infty \le \beta$ allows us to avoid issues of ever-increasing search directions. Overall, in the sequel we denote by $\widetilde {\mathcal{X}}(x)$ and $d(x)$ the convex feasible set and the unique solution of subproblem \eqref{eq:p_k}, respectively, i.e.
\[
\begin{array}{rcl}
\widetilde {\mathcal{X}}(x) & \triangleq & \left\{d \in \mathbb R^n \, : \, \tilde g(d; x) \le \kappa(x) \, e, \, \|d\|_\infty \le \beta, \, d \in K - x \right\}\\[5pt]
d(x) & \triangleq & \argmin_d \{\tilde f(d; x) \, | \, d \in \widetilde{\mathcal{X}}(x)\},
\end{array}
\]
and we equivalently refer to constraint $\|d\|_\infty \le \beta$ as $d \in \beta \mathbb B^n_\infty$, where $\mathbb B^n_\infty$ is the closed unit ball in $\mathbb R^n$ associated with the infinity-norm.

For our approach to be legitimate and lead to useful convergence results, we obviously need to make assumptions on the surrogate functions $\tilde f$ and $\tilde g$.

\medskip \noindent {\bf{Assumption A}}

\medskip

\noindent {\textit{Let $O_d$ and $O_x$ be open neighborhoods of $\beta \mathbb B^n_\infty$ and $K$, respectively, and $\tilde{f}:O_d\times O_x\rightarrow\mathbb{R}$ and $\tilde g_i:\mathbb R^n\times O_x\rightarrow\mathbb R$ for every $i = 1, \ldots, m$ be continuously differentiable on $O_d$ with respect to the first argument and such that}}
\begin{description}
\item [{\textit{A1)}}] {\textit{$\tilde{f}(\bullet;x)$ is a strongly convex function on $O_d$ for every $x \in K$ with modulus of strong convexity $c >0$ independent of $x$;}}

\item [{\textit{A2)}}] {\textit{$\tilde f(\bullet;\bullet)$ is continuous on $O_d \times O_x$;}}

\item [{\textit{A3)}}] {\textit{$\nabla_1 \tilde{f}(\bullet;\bullet)$ is continuous $O_d \times O_x$;}}

\item [{\textit{A4)}}] {\textit{$\nabla_1 \tilde f(0; x) = \nabla f(x)$ for every $x \in K$;}}

\item [{\textit{A5)}}] {\textit{$\tilde g_i(\bullet;x)$ is a convex function on $O_d$ for every $x \in K$;}}

\item [{\textit{A6)}}] {\textit{$\tilde g_i(\bullet;\bullet)$ is continuous on $\mathbb R^n \times O_x$;}}

\item [{\textit{A7)}}] {\textit{$\tilde g_i(0;x) = g_i(x)$ for every $x \in K$;}}

\item [{\textit{A8)}}] {\textit{$\nabla_1 \tilde g_i(\bullet;\bullet)$ is continuous on $O_d \times O_x$;}}

\item [{\textit{A9)}}] {\textit{$\nabla_1 \tilde g_i(0;x)=\nabla g_i(x),$ for every $x \in K$;}}
\end{description}
{\textit{where $\nabla_1 \tilde f(u;x)$ and $\nabla_1 \tilde g_i(u;x)$ denote the partial gradient of $\tilde f(\bullet;x)$ and $\tilde g_i(\bullet;x)$ evaluated at $u$.}}
These conditions are easily satisfied in practice and have been employed in many recent papers. While we refer the reader to \cite{facchinei2017feasible,scutari2014decomposition} as good sources of examples, we nevertheless pause to consider some possible choices for the surrogate functions $\tilde f$ and $\tilde g$ and to make a few general considerations.

\subsection{On the choice of $\tilde f$ and $\tilde g$}\label{sec:subsec3.1}
The direction finding subproblem \eqref{eq:p_k} is a direct generalization of traditional SQP subproblems and, in particular, of the subproblems considered in \cite{burke1989sequential}. The most classical choice for $\tilde f$ and $\tilde g$ are
\begin{equation}\label{eq:quad_lin_approx}
\tilde f (d; x) \triangleq \nabla f(x)\trt d + \frac{1}{2} d\trt H(x) d, \qquad \quad
\tilde g (d; x) \triangleq g(x) + \nabla g(x)\trt d,
\end{equation}
where $H(x)$ is a positive definite symmetric matrix. With this choice, Assumption A can be easily  satisfied provided that,  classically,  the smallest eigenvalue of the positive definite matrix $H(x)$ is uniformly bounded away from zero. Note that if we use the surrogates \eqref{eq:quad_lin_approx} in \eqref{eq:p_k}, and assuming $K = \mathbb R^n$, \eqref{eq:p_k} becomes the more classical SQP-type subproblem
\[
\begin{array}{cl}
\underset{d}{\mbox{minimize}} & f(x) + \nabla f(x)\trt d + \frac{1}{2} d\trt H(x) d\\
\mbox{s.t.} &  g(x) + \nabla g(x)\trt d \le \kappa(x)\\[5pt]
& \|d\|_\infty \le \beta.\\[5pt]
\end{array}
\]
Regarding $\kappa$, we remark that, with the classical choice given in \eqref{eq:quad_lin_approx}, problem \eqref{eq:capoptpro} in its definition reduces to a linear program if $K$ is polyhedral, and thus can be efficiently  solved.

While the discussion above shows that we can cover classical SQP schemes by using linear/quadratic approximations, it is interesting to at least hint at how the flexibility allowed by Assumption A can be exploited to define better approximations to the  original problem \eqref{eq:startpro}. Suppose for example that $f(x) = f_1(x) + f_2(x)$ with both functions $C^{1,1}$, but with $f_1$ convex and $f_2$ not necessarily so. Instead of approximating the whole function with a quadratic model, we could well ``preserve'' the convex part  and only approximate the nonconvex one, therefore setting
\[
\tilde f(d;x) = f_1(x+d) + f_2(x) +  \nabla f_2(x)\trt d + \frac{1}{2} d\trt H(x) d
\]
with $H(x)$ as before. It is clear that this $\tilde f$ satisfies Assumption A and is presumably a better approximation to $f$ than $ \nabla f(x)\trt d + \frac{1}{2} d\trt H(x) d$ considered above.

As a further example, assume that $f$ is the product of two functions $f_1(x)f_2(x)$ with $f_1$ and $f_2$ convex and (for simplicity of presentation) positive. This is a rather frequent case in applications, see for example  \cite{scutari2017parallelbis}. Since we have $\nabla f(x) = f_2(x) \nabla f_1(x) + f_1(x) \nabla f_2(x)$, it seems rather natural to set
\[
\tilde f(d;x) = f_2(x) f_1(x+d) + f_1(x)  f_2(x+d) + \frac{1}{2} d\trt H(x) d,
\]
which, again, should result in a sensibly tailored approximation that preserves part of the structure of the objective function.

Of course, an underlying assumption of our approach is that  subproblem \eqref{eq:p_k} can be solved
efficiently. We do not insist on this point because it is very dependent on the choice of $\tilde f$ and $\tilde g$,
which in turn is dictated by the original problem \eqref{eq:startpro}. But the use of models that go beyond the
classical quadratic/linear one in constrained optimization is emerging consistently in the literature since it
permits one to exploit any potentially favorable structure in problem \eqref{eq:startpro} and, in any case, to better
tailor the subproblems to the original problem, see for example
\cite{beck2010sequential,facchinei2017feasible,lipp2016variations,scutari2014decomposition,
sun2017majorization,svanberg2002class}. This use is also motivated by the
possibility to solve efficiently more complex subproblems than the classical quadratic ones, sometimes even in
closed form, and by the desire
for faster convergence behaviors, see e.g. \cite{facchinei2017feasible,hong2016unified,mairal2013optimization,martinez2017high,
scutari2017parallelbis,razaviyayn2013unified,sun2017majorization} and references therein.
\smallskip

Among all the possible choices for the approximating functions, the case where we take the $\tilde g_i$s to be
Upper Convex Approximations (UCA) of $g_i$s is worth to be pointed out. More precisely, suppose that, in
addition to Assumption A, for every $x \in K$ we have
\begin{equation}\label{eq:upperb}
\tilde g_i(d;x)\, \geq \,g_i(x+d), \;\; \forall i=1,\ldots, m, \;\; \forall d\,\in\, K-x.
\end{equation}
The main consequence  of this choice is that if $x \in \mathcal X$, then $0 \in \widetilde{\mathcal X}(x)$ and,
by \eqref{eq:upperb}, $x + \widetilde {\mathcal X}(x) \subseteq \mathcal X$.
This means that if $x^\nu \in {\cal X}$ and, according to \eqref{eq:method}, we set $x^{\nu+1} = x^\nu +
\gamma^\nu d(x^\nu)$  with $\gamma^\nu \in (0, 1]$ (a condition that will always be satisfied by all
algorithms considered in this paper), also $x^{\nu +1}$ is feasible, i.e.  $x^{\nu +1}\in {\cal X}$.
This simple observations has important algorithmic ramifications that will be explored further in the next three sections.
The main issue if one wants to use UCAs is finding suitable majorants. It turns out
this can be done in a host of situations; the interested reader can find a
very rich array of examples in
\cite{facchinei2017feasible,scutari2017parallel,scutari2017parallelbis,hong2016unified,hunter2004tutorial,
razaviyayn2013unified,sun2017majorization}.
Here,   we just consider two examples.

The simplest case is possibly the one in
\cite{auslender2010moving}, where a feasible SQP-like approach is developed that rests on the assumption
(among others) that the $g_i$  have Lipschitz continuous gradients (with constant $L_{\nabla g_i}$) and the
descent lemma is
used for defining suitable majorants as
 \begin{equation}\label{eq:UCAg}
 \tilde g_i(d; x) = g_i(x) + \nabla g_i(x)\trt d + \frac{a}{2}\| d\|^2.
 \end{equation}
By taking $a \ge L_{\nabla g_i}$, \eqref{eq:UCAg} provides an UCA for $g_i$.
A second example of a case in which we can very easily build majorants is when the function has a DC structure.
Specifically, suppose  that  $g_i = g_i^+ - g_i^-$, with both $g_i^+$ and $ g_i^-$ convex. In this case we can
build an upper convex approximation by setting
\[
\tilde g_i(d;x) = g_i^+(x+d) - (g_i^-(x) + \nabla g_i^-(x)\trt d).
\]

\subsection{Main properties of \eqref{eq:p_k}} 
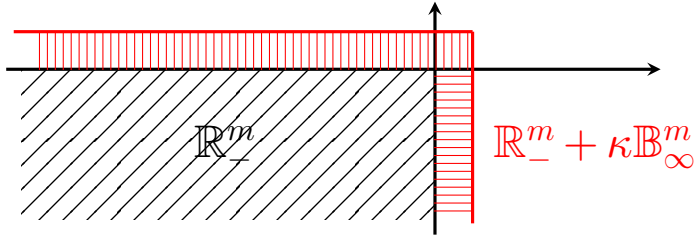
\begin{figure}[h]
\centering
\begin{tikzpicture}

\tikzset{
hatch distance/.store in=\hatchdistance,
hatch distance=10pt,
hatch thickness/.store in=\hatchthickness,
hatch thickness=2pt
}

\makeatletter
\pgfdeclarepatternformonly[\hatchdistance,\hatchthickness]{flexible hatch}
{\pgfqpoint{0pt}{0pt}}
{\pgfqpoint{\hatchdistance}{\hatchdistance}}
{\pgfpoint{\hatchdistance-1pt}{\hatchdistance-1pt}}%
{
\pgfsetcolor{\tikz@pattern@color}
\pgfsetlinewidth{\hatchthickness}
\pgfpathmoveto{\pgfqpoint{0pt}{0pt}}
\pgfpathlineto{\pgfqpoint{\hatchdistance}{\hatchdistance}}
\pgfusepath{stroke}
}
\makeatother

\coordinate (origo) at (0,0);
\coordinate (pivot1) at (0.5,0.5);
\coordinate (pivot2) at (0.5,0);

\draw[black,very thick,->,>=stealth] (origo) -- ++(3,0) node[black,right] {};
\draw[black,very thick,-] (origo) -- ++(-5.7,0) node[black,right] {};
\draw[black,very thick,->,>=stealth] (origo) -- ++(0,0.9) node[black,above] {};
\draw[black,very thick,-] (origo) -- ++(0,-2.2) node[black,above] {};

\draw[red,very thick,-] (pivot1) -- ++(-6.1,0) node[black,right] {};
\draw[red,very thick,-] (pivot1) -- ++(0,-2.55) node[black,right] {};


\fill[pattern=flexible hatch,
hatch distance=10pt,
hatch thickness=0.5pt] ($ (origo) + (-5.5,-2) $) rectangle ($ (origo)$);
\fill[pattern = vertical lines,pattern color=red] ($(-5.3,0) $) rectangle ($(pivot1)$);
\fill[pattern = horizontal lines,pattern color=red] ($(0,-1.9) $) rectangle ($(0.5,0)$);

\node[text width=40pt,] at (-2.5,-1.05) {\Large{$\mathbb R^{m}_{-}$}};

\node[text width=115pt,] at (2.8,-1.05) {\textcolor{red}{\Large{$\mathbb R^{m}_{-} + \kappa \mathbb B^{m}_{\infty}$}}};

\end{tikzpicture}
\caption{From $\mathbb R_-^m$ to $\mathbb R_-^m + \kappa \mathbb B^m_\infty$: the enlargement in the feasible region of \eqref{eq:p_k}\label{fig:collar}}
\end{figure}

In this section we state the main properties of subproblem \eqref{eq:p_k}. Starting from feasibility, we remark, as already mentioned, that the term $\kappa(x)$ plays a key role in guaranteeing that our subproblems \eqref{eq:p_k} have a nonempty feasible set $\widetilde {\mathcal{X}}(x)$. Since $\kappa (x)$ is always nonnegative, being the sum of two nonnegative quantities, it restores feasibility by enlarging (with respect to the SQP choice $\tilde g(d;x) \leq 0$) the range of admissible values, see Figure
\ref{fig:collar}.
In fact, the feasible set of problem \eqref{eq:p_k}, for every $x \in K$, is nonempty: choosing $\hat d$ at which the minimum in the expression of $\kappa(x)$ is attained, we have
$$
\tilde g(\hat d; x) \le \min_d \left\{\max_i\{\tilde g_i(d; x)_+\} \, | \, \|d\|_\infty \le \rho, \, d \in K - x\right\} e = \max_i\{\tilde g_i(\hat d; x)_+\} e,
$$
and, in turn,
$$
\begin{array}{rcl}
\tilde g(\hat d; x) & = & (1 - \lambda) \tilde g(\hat d; x) + \lambda \tilde g(\hat d; x)\\[5pt]
& \le & \left[(1 - \lambda)  \max_i\{\tilde g_i(0; x)_+\} + \lambda \min_d \left\{\max_i\{\tilde g_i(d; x)_+\} \, | \, \|d\|_\infty \le \rho, d \in K - x\right\}\right] e = \kappa(x) e.
\end{array}
$$
In Lemma \ref{th:prelresfeas}, we establish some preliminary properties concerning the feasible set of problem \eqref{eq:p_k}.
\begin{lemma}\label{th:prelresfeas}
\begin{enumerate}
\item[(i)]
For every $\hat x \in K$, and for every $\alpha > 0$ and $d \in \alpha \mathbb B_\infty^n \cap (K - \hat x)$, the constraint qualification
\begin{equation}\label{eq:furcq}
[-N_{\alpha \mathbb B^n_\infty}(d)] \cap N_{K-\hat x}(d) = \{0\}
\end{equation}
holds and, in turn, $N_{\alpha \mathbb B^n_\infty \cap (K - \hat x)} (d) = N_{\alpha \mathbb B_\infty^n}(d) + N_{K - \hat x}(d)$;
\item[(ii)]
for every $\alpha > 0$, the set-valued mapping $\alpha \mathbb B_\infty^n \cap (K - \bullet)$ is continuous on  $K$ relative to $K$;
\item[(iii)]
letting $C \triangleq \{(d, x) \in \beta \mathbb B_\infty^n \times K \, : \, d + x \in K \}$, the set-valued mapping $N_{\beta \mathbb B_\infty^n \cap (K - \bullet)}(\bullet)$ is outer semicontinuous on $C$ relative to $C$.
\end{enumerate}
\end{lemma}
\proof{Proof.}
(i) Let $0 \neq \eta \in [- N_{\alpha \mathbb B_\infty^n}(d)] \cap N_{K- \hat x}(d)$. Thanks to the convexity of the sets $\alpha \mathbb B_\infty^n$ and $K - \hat x$, we have $- \eta\trt (v - d) \le 0$ $\forall v \in \alpha \mathbb B^n_\infty$ and $\eta\trt (y - d) \le 0$ $\forall y \in (K - \hat x)$. Choosing $y = 0 \in (K - \hat x)$, one gets the following contradiction:
\[
0 < \alpha \, \max_v \{-\eta\trt v \, | \, v \in \mathbb B^n_\infty\} \le -\eta\trt d \le 0,
\]
thus proving relation \eqref{eq:furcq}. As a consequence, the other claim in (i) follows from \cite[Theorem 6.42]{RockWets98}.

(ii)  The property is due to the continuity (relative to $K$) of the set-valued mapping $K - \bullet$ at every $x \in K$ and to the fact that $\alpha \mathbb B_\infty^n \cap (K - x) \neq \emptyset$ for every $x \in K$.

(iii) Suppose by contradiction that $(d^\nu, x^\nu) \underset{C}{\to} (\bar d, \bar x)$, $\ \eta^\nu \in N_{\beta \mathbb B_\infty^n \cap (K - x^\nu)}(d^\nu)$, $\eta^\nu \to \bar \eta$ with $\bar \eta \notin N_{\beta \mathbb B_\infty^n \cap (K - \bar x)}(\bar d)$. Hence, $\bar z \in \beta \mathbb B_\infty^n \cap (K - \bar x)$ exists such that $\bar \eta\trt (\bar z - \bar d) > 0$. By the inner semicontinuity relative to $K$ (see \cite{RockWets98} for the definition of inner semicontinuity) of $\beta \mathbb B_\infty^n \cap (K - \bullet)$ at $\bar x$, $z^\nu$ exists such that $z^\nu \to \bar z$ and $z^\nu \in \beta \mathbb B_\infty^n \cap (K - x^\nu)$. In turn, eventually we get $(\eta^\nu)\trt (z^\nu - d^\nu) > 0$ in contradiction to the inclusion $\ \eta^\nu \in N_{\beta \mathbb B_\infty^n \cap (K - x^\nu)}(d^\nu)$.
\hfill  \Halmos\endproof
The function $\kappa(x)$ is obviously continuous and, under a very weak additional requirement, also locally Lipschitz continuous. This result has been shown in \cite{burke1989sequential} whenever
$\tilde g$ is the linear approximation in \eqref{eq:quad_lin_approx} and readily generalizes to the case of the surrogate $\tilde g$ we consider here.\begin{proposition}\label{eq:kappa lip}
Under Assumption A, $\kappa(\bullet)$ is  continuous on $K$ relative to $K$. If, in addition, $\tilde g(\bullet; \bullet)$ is locally Lipschitz continuous on $O_d \times O_x$, then $\kappa(\bullet)$ is also locally Lipschitz continuous on $K$.
\end{proposition}
\proof{Proof.}
The continuity of  $\kappa(\bullet)$ follows readily from the continuity (relative to $K$) of the set-valued mapping $\rho \mathbb B_\infty^n \cap (K - \bullet)$ at every $x \in K$: this in turn follows from (ii) in Lemma \ref{th:prelresfeas} with $\alpha = \rho$.

The Lipschitz continuity under the additional condition derives from, e.g., \cite[Theorem 3.1]{Rock85}. Suffice it to observe that the constraint qualification \eqref{eq:furcq} with $\alpha = \rho$ holds for every $x \in K$ and $d \in \rho \mathbb B_\infty^n \cap K - x$, and the problem in the definition of $\kappa$ is solvable for every $x$ in a neighborhood of $K$.
\hfill  \Halmos\endproof

Note that the local Lipschitz continuity of $\tilde g(\bullet; \bullet)$ is part of Assumption C to be introduced shortly.

The following technical lemma is very useful for the subsequent developments.
\begin{lemma}\label{th:prellem}
Under Assumption A, the following results hold for any $\hat x \in K$:
\begin{enumerate}
\item[(i)]
if $\max_i \{g_i(\hat x)_+\} > 0$ and $\kappa(\hat x) < \max_i \{g_i(\hat x)_+\}$, then, for all $\rho \in (0, \beta)$, there exists $d \in \mathrm{int} (\beta \mathbb B^n_\infty) \cap \mathrm{rel \, int}(K - \hat x)$ such that $\tilde g(d; \hat x) < \kappa(\hat x) e$;

\item[(ii)]
if $\max_i \{g_i(\hat x)_+\} > 0$ and $\kappa(\hat x) = \max_i \{g_i(\hat x)_+\}$, then $\hat x$ is an ES point for \eqref{eq:startpro};

\item[(iii)]
if $\max_i \{g_i(\hat x)_+\} = 0$, then either $\hat x$ is a FJ point for \eqref{eq:startpro} or, for all $\rho \in (0, \beta)$, there exists $d \in \mathrm{int} (\beta \mathbb B^n_\infty) \cap \mathrm{rel \, int}(K - \hat x)$ such that $\tilde g(d; \hat x) < 0$.
\end{enumerate}
\end{lemma}
\proof{Proof.}
(i) Choosing $\hat d \in \argmin_d \left\{\max_i\{\tilde g_i(d; \hat x)_+\} \, | \, \|d\|_\infty \le \rho, \, d \in K - \hat x \right\}$, we can infer $\tilde g(\hat d; \hat x) \le \min_d \left\{\max_i\{\tilde g_i(d; \hat x)_+\} \, | \, \|d\|_\infty \le \rho, d \in K - \hat x \right\} e$, while $\tilde g(\hat d; \hat x) \le \kappa(\hat x) e < \max_i\{g_i(\hat x)_+\} e$ and, thus,
$$
\tilde g(\hat d; \hat x) = \lambda \tilde g(\hat d; \hat x) + (1 - \lambda) \tilde g(\hat d; \hat x) < \kappa(\hat x) e,
$$
with $\hat d \in \rho \mathbb B_\infty^n \cap (K - \hat x)$. The claim follows by continuity since $\rho < \beta$.

(ii) Equality $\kappa(\hat x) = \max_i\{g_i(\hat x)_+\}$ holds if and only if $d = 0$ solves the minimization problem in the definition of $\kappa$ and, in turn, $M_0(\hat x) \neq \{0\}$ by \eqref{eq:furcq} with $\alpha = \rho$, A7 and A9.

(iii) With $\max_i \{g_i(\hat x)_+\}$ being equal to zero, we have $\kappa(\hat x) = 0$ and $g(\hat x) \le 0$. If $M_0(\hat x) \neq \{0\}$, then, by definition, $\hat x$ is a FJ point for \eqref{eq:startpro} and the result holds.

Thus, let us suppose $M_0(\hat x) = \{0\}$. For those $j \in \{1, \ldots, m\}$ such that $g_j(\hat x) < 0$, we have $\tilde g_j(0; \hat x) = g_j(\hat x) < 0$; as for indices $k \in \{1, \ldots, m\}$ with $g_k(\hat x) = 0$, by \eqref{eq:EMFCQ1}, there exists $\hat d \in T_K(\hat x)$ such that
$$
0 > \nabla g_k(\hat x)\trt \hat d = \nabla_1 \tilde g_k(0; \hat x)\trt \hat d =  \lim_{\tau \downarrow 0} \frac{\tilde g_k(\tau \hat d; \hat x) - \tilde g_k(0; \hat x)}{\tau}.
$$
Thus, there exists a sequence $\{d^\nu\}$ of feasible directions for $K$ at $\hat x$ such that $d^\nu \in T_K(\hat x)$ and $d^\nu \to \hat d$. Choosing $\tau^\nu$ sufficiently small, we get $\hat x + \tau^\nu d^\nu \in K$ for every $\nu$ and the claim follows by continuity, observing that $\tilde g_i(\tau \hat d; \hat x) < 0$ for every $i$ and for any $\tau$ sufficiently small.
\hfill  \Halmos\endproof
The quantity
\begin{equation}\label{eq:theta}
\theta(x) \triangleq  \max_i \{g_i(x)_+\} - \kappa(x) = \lambda \left(\max_i \{g_i(x)_+\} - \min_d \left\{\max_i  \left\{\tilde g_i(d; x)_+ \right\} \, | \, \|d\|_\infty \le \rho, \, d \in K - x\right\}\right)
\end{equation}
plays a key role in the previous lemma and in all the subsequent developments. As shown in the following proposition, $\theta$ turns out to be a stationarity measure for the violation-of-the-constraints problem \eqref{eq:feaspro}.
\begin{proposition}\label{th:theta}
Under Assumption A,
\begin{itemize}
\item[(i)]	
the nonnegative function $\theta(\bullet)$ is continuous on $K$ relative to $K$;
\item[(ii)]
$\theta(\hat x) = 0$ if and only if $\hat x$ is a stationary point for problem \eqref{eq:feaspro};
\item[(iii)]
we have, for every $x \in K$,
\begin{equation}\label{eq:thetadelta}
\theta(x) \le \|\nabla g(x)\|_\infty \|d(x)\|.
\end{equation}
\end{itemize}
\end{proposition}
\proof{Proof.}
(i) By the definition \eqref{eq:cap} of $\kappa$, $\kappa(\hat x) \le \max_i\{g_i(\hat x)_+\}$ since $d = 0$ is feasible for the minimization problem in \eqref{eq:cap} and A7 holds. The continuity follows from Proposition \ref{eq:kappa lip}.

(ii) It is also clear that at any feasible point $\hat x$ for \eqref{eq:startpro}, $\theta(\hat x) = 0$; of course, every feasible point for \eqref{eq:startpro} is stationary for problem \eqref{eq:feaspro}. Consider now an infeasible point $\hat x$ for \eqref{eq:startpro} and suppose that $\theta(\hat x)=0$. By (ii) in Lemma \ref{th:prellem}, $\hat x$ turns out to be an ES point for \eqref{eq:startpro}. Hence, we are left to show that if $\hat x$ is an ES point for \eqref{eq:startpro}, then $\theta(\hat x) = 0$. For $\hat x$ to be ES, it is necessary and sufficient (see condition \eqref{eq:statfeaspro}) to have $M_0(\hat x) \neq \{0\}$ which in turn, by the Motzkin's alternative theorem (see e.g. \cite[2.5.2]{craven2012mathematical}), holds if and only if
\begin{equation}\label{eq:Motzk}
\nexists \, d \in T_K(\hat x) \, : \, \nabla g_i(\hat x)\trt d < 0, \, \forall i \in I_+(\hat x) \triangleq \{i \, : \, g_i(\hat x) = \max_j \{g_j(\hat x)_+\}\}.
\end{equation}
Suppose by contradiction that $\theta(\hat x) > 0$. Then, noting that $d\in K-\hat x$ implies $d\in T_K(\hat x)$, Lemma \ref{th:prellem} (i) states that $d \in T_K(x)$ exists such that
$\tilde g_i(d;\hat x) < \kappa(\hat x)$ for all $i \in I_+(\hat x)$. But then, using A5, A7, and A9, we can write, for every $i \in I_+(\hat x)$,
\[
\max_i \{g_i(\hat x)_+\}  \,>\, \kappa(\hat x) \,> \,\tilde g_i(d;\hat x)  \,\geq\, \tilde g_i(0;\hat x) + \nabla_1 \tilde g_i(0;\hat x)\trt(d-0) \,\geq\,  g_i(\hat x) + \nabla g_i(\hat x)\trt d.
\]
Since $i\in I_+(\hat x)$, this implies  $\nabla g_i(\hat x)\trt d <0$, contradicting \eqref{eq:Motzk}.

(iii) Furthermore,
\begin{equation*}\label{eq:phicappa}
\begin{array}{rcl}
0 & \le & \theta(x^\nu) = \displaystyle \max_i\{g_i(x^\nu)_+\} - \kappa(x^\nu) \overset{(a)}{\le} \displaystyle \max_i\{g_i(x^\nu)_+\} - \max_i\{\tilde g_i(d(x^\nu); x^\nu)_+\}\\[5pt]
& \overset{(b)}{\le} & \displaystyle \max_i\{g_i(x^\nu)_+\} - \max_i\{(g_i(x^\nu) +
 \nabla g_i(x^\nu)\trt d(x^\nu))_+\}\\[5pt]
& \overset{(c)}{\le} & \displaystyle \max_i\{(g_i(x^\nu) - g_i(x^\nu) - \nabla g_i(x^\nu)\trt d(x^\nu))_+\} \le \|\nabla g(x^\nu)\trt d(x^\nu)\|_\infty \le \|\nabla g(x^\nu)\|_\infty \|d(x^\nu)\|,
\end{array}
\end{equation*}
where (a) holds since $\tilde g(d(x^\nu); x^\nu) \le \kappa(x^\nu) e$, (b) is due to A5, A7 and A9, and (c) follows observing that $\max\{0, \alpha_1\} - \max\{0, \alpha_2\} \le \max\{0, \alpha_1 - \alpha_2\} $ for any $\alpha_1, \, \alpha_2 \in \mathbb R$.
\hfill  \Halmos\endproof

Leveraging Lemma \ref{th:prellem} and resorting to standard results in parametric optimization, we can establish a key continuity property for the solution mapping $d(\bullet)$ of subproblem \eqref{eq:p_k}.

\begin{proposition}\label{th:keyemfcq}
Under Assumption A, let the eMFCQ hold at $\hat x \in K$ for problem \eqref{eq:startpro}. Then,

(i) the MFCQ holds at every point of $\widetilde {\mathcal{X}}(\hat x)$ for subproblem {\emph{(P$_{\hat x}$)}};

(ii)  a neighborhood $\mathcal V$ of $\hat x$ exists such that, for every point $x \in K \cap \mathcal V$, the function $d(\bullet)$ is continuous relative to $K$.
\end{proposition}
\proof{Proof.}
If the eMFCQ holds at $\hat x$ for problem \eqref{eq:startpro}, case (ii) in Lemma \ref{th:prellem} cannot occur. On the other hand, as for both cases
(i) and (iii) in Lemma \ref{th:prellem}, Slater's constraint qualification holds for $\widetilde {\mathcal{X}}(\hat x)$ and, since
$\widetilde {\mathcal{X}}(\hat x)$ is convex, this proves (i).
Thanks to A6 and (ii) in Lemma \ref{th:prelresfeas}, the set-valued mapping $\widetilde
{\mathcal{X}}(\bullet) = [\beta \mathbb B^n_\infty \cap (K - \bullet)] \cap \left\{d \in \mathbb R^n \, : \, \tilde g(d; \bullet) \le \kappa(\bullet)
e\right\}$ is outer semicontinuous at $\hat x$ relative to $K$ by \cite[Theorem 3.1.1]{banknon}, having taken into account  that  $
\kappa(\bullet)$ is continuous  by Proposition \ref{eq:kappa lip}. Moreover, $\widetilde {\mathcal{X}}(\bullet)$, by
virtue of the Slater's constraint qualification, A5, A6 and (ii) in Lemma \ref{th:prelresfeas}, is
also inner semicontinuous (see \cite[Theorem 3.1.6]{banknon}) at $\hat x$ relative to $K$. Hence, thanks to A1, the continuity (relative to $K$) of
$d(\bullet)$, leveraging \cite[Theorems 3.1.1 and 4.3.3]{banknon}, follows from \cite[Corollary 5.20]{RockWets98}.
\hfill  \Halmos\endproof

To prove  some refinements of the convergence results in the next section,  we need
$d(\bullet)$ to be not only continuous, but also H\"older continuous on compact sets: for this reason, we introduce Assumption B.

\medskip \noindent {\bf{Assumption B}}

\medskip

\noindent {\textit{For any compact set $S\subseteq K$, two positive constants $\theta$ and $\alpha$ exist such that}}
$$
\| d(y) - d(z)\| \leq \theta \|y-z\|^\alpha, \quad \forall y, z \in S.
$$
Since it is not immediately obvious when this condition is satisfied, below we give a set of simple sufficient conditions on
$\tilde f$ and $\tilde g$ for Assumption B to hold.

\medskip \noindent {\bf{Assumption C}}

\medskip

{\begin{description}
\item [{\textit{C1)}}] {\textit{$\nabla_1 \tilde f(\bullet;\bullet)$ is locally Lipschitz continuous on $O_d \times O_x$;}}

\item [{\textit{C2)}}] {\textit{each $\tilde{g}_j(\bullet;\bullet)$ is locally Lipschitz continuous on $O_d \times O_x$.}}
\end{description}}

The following proposition, which builds on the results in \cite{yen95}, shows the desired result.
\begin{proposition}\label{th:uspropofupgpre}
Under Assumptions A and C, let $S\subseteq K$ be compact. Suppose further that the eMFCQ
holds
at every $\hat x \in S$. Then, there exists $\theta > 0$ such that, for every ${y},{z}\in S$,
\begin{equation}\label{eq:hol}
\|d({y})-d({z})\|\le\theta\|{y}-{z}\|^{\frac{1}{2}}.
\end{equation}

\end{proposition}

\proof{Proof.}
Preliminarily, observe  that by Proposition \ref{eq:kappa lip}, $\kappa(\bullet)$ is locally Lipschitz continuous. Furthermore, by Proposition \ref{th:keyemfcq} (i), we have that the MFCQ holds at every point in $\widetilde{\cal X}(\hat x)$ and, in particular at $d(\hat x)$. In turn,
the MFCQ at $d(\hat x) \in \widetilde{\mathcal X}(\hat x)$, for every $\hat x \in S$, implies, by \cite[Theorem 3.2]{Rock85},
that the set-valued mapping $\widetilde{\mathcal X}$ has the Aubin property relative to $S$
at $\hat x$ for $d(\hat x)$
for every $\hat x \in S$ (see \cite{RockWets98} for the definition of the Aubin property). Therefore, in view of \cite[Theorem 2.1]{yen95}, for every $\hat x \in S$, there exist $\theta' > 0$, $\theta'' > 0$ and a neighborhood ${\mathcal{V}}$ of $\hat x$ such that,
for every $y, \, z \in {\mathcal{V}}\cap S$
\begin{equation*}
\|d({y})-d({z})\|\le\theta' \|{y}-{z}\|+\theta'' \|{y}-{z}\|^{\frac{1}{2}}.
\end{equation*}
By the previous relation and the compactness of set $S$, \eqref{eq:hol} holds. \hfill
\hfill  \Halmos\endproof

\begin{remark}
Assumptions A and C may look tediously detailed, but this is necessary to correctly identify the minimal conditions that make our method work.
We emphasize that these conditions are trivially satisfied when one uses as $\tilde f$ and $\tilde g$ the classical quadratic/linear approximations \eqref{eq:quad_lin_approx} of standard SQP methods.
Assumption C reinforces some of the requirements in Assumption A; we refer the reader to
\cite{facchinei2017feasible} for some examples of surrogate $\tilde g$s satisfying (Assumption A and)
Assumption C beyond the obvious case of linear approximations.
\end{remark}

We conclude this section discussing the KKT conditions for problem \eqref{eq:p_k}.
Observe preliminarily that the constraint $\| d\|_\infty \leq \beta$ corresponds to $2n$ bounds  of the type
$-\beta \leq d_i \leq \beta$. However, in what follows we are interested only in the multipliers corresponding to the constraints $\tilde g(d;x) \leq \kappa(x) e$, and therefore we find it expedient to write the KKT conditions as
\begin{equation*}\label{eq:kktindicator}
0 \in \nabla_1 \tilde f(d(x); x) + \nabla_1 \tilde g(d(x); x)  \xi + N_{\beta \mathbb B^n_\infty \cap (K - x)}(d(x)),
\end{equation*}
with the KKT multipliers $\xi$ satisfying the conditions $\xi \geq 0$ and $\xi\trt (\tilde g(d(x); x)- \kappa(x) e) =0$. We now establish the local boundedness of these KKT multipliers.

\begin{proposition}\label{lem:boundedness}
Under Assumption A, let $\hat x$ belong to $K$ and suppose that $\hat d \in \beta \mathbb B^n_\infty \cap (K - \hat x)$ exists such that $\tilde g(\hat d; \hat x) < \kappa(\hat x) e$.
Then, a neighborhood $\mathcal V$ of $\hat x$ exists such that, for every point $x \in K \cap \mathcal V$, the unique solution $d(x)$ of \eqref{eq:p_k} is a KKT point for problem \eqref{eq:p_k} and the set-valued mapping of the KKT multipliers is locally bounded at $\hat x$ relative to $K$.
\end{proposition}

\proof{Proof.}
The condition  $\tilde g(\hat d; \hat x) <\kappa(\hat x) e$ with $\hat d \in \beta \mathbb B^n_\infty \cap (K - \hat x)$ is nothing else but the Slater's CQ for problem (P$_{\hat x}$), which obviously implies that the MFCQ holds at the unique solution of problem (P$_{\hat x}$). The derivation of the result is then rather classical and follows easily from, e.g., \cite[Proposition 5.4.3]{FacchPangBk} taking into account Lemma \ref{th:prelresfeas} (ii),  Propositions \ref{eq:kappa lip} and \ref{th:keyemfcq}, and the outer semicontinuity of $N_{\beta \mathbb B_\infty^n \cap (K - \bullet)}(\bullet)$, see Lemma \ref{th:prelresfeas} (iii).
\hfill  \Halmos\endproof

\section{Convergence of DSMs}\label{Sec:convergence}
We are now ready to introduce the proposed scheme, as given in Algorithm \ref{algoBasic}.
{\setlength{\skiptext}{10pt}
\setlength{\skiprule}{5pt}

\IncMargin{1em}
\begin{algorithm}
\KwData{$\gamma^\nu \in (0,1]$ such that \eqref{eq:gamma} holds, $x^{0} \in K$,  $\nu \longleftarrow 0$\;}
\Repeat{
{\nlset{(S.1)} \If{$x^{\nu}$ {\emph{is generalized stationary for}} \eqref{eq:startpro}}{{\bf{stop}} and {\bf{return}} $x^\nu$\;} \label{S.11}}
\nlset{(S.2)} compute $\kappa(x^\nu)$ and the solution $d(x^\nu)$ of problem (P$_{x^\nu}$)\; \label{S.12}
\nlset{(S.3)} set $x^{\nu+1}=x^{\nu}+\gamma^{\nu}d(x^\nu)$, $\nu\longleftarrow\nu+1$\; \label{S.13}}{}
\caption{\label{algoBasic} DSM Algorithm for \eqref{eq:startpro}}
\end{algorithm}}
\DecMargin{1em}
\noindent The algorithm is always well defined if Assumption A, which guarantees existence and uniqueness of $d(x^\nu)$,  holds.
The main (and essentially only) computational burden is
given by the computation of $\kappa(x^\nu)$ and the solution of the strongly convex subproblem (P$_{x^\nu}$). This difficulty can range from that necessary to solve an LP and a strongly convex quadratic problem, whenever quadratic/linear approximations are used, to that of solving two convex optimization problems. Theorem \ref{th:converg} below establishes the main convergence properties of Algorithm \ref{algoBasic}.
In a nutshell, the theorem shows that, unless $x^\nu$ is a generalized stationary point, $d(x^\nu)$ is a descent direction for $W(x^\nu; \varepsilon)$ if $\varepsilon$ is sufficiently small. Elaborating on this simple fact we can then show, without ever computing $W$ or actually determining a value for $\varepsilon$, that the sequence generated eventually lands on a generalized stationary point.
The results in Theorem \ref{th:converg} do not exclude the   possibility that Algorithm \ref{algoBasic} generates an unbounded sequence. In Section \ref{sec:bounded} we discuss the meaning of this possible outcome and, more importantly, give several conditions under which we can guarantee that the sequence generated by Algorithm \ref{algoBasic} (and also by the two algorithms we introduce in the next section) is bounded.

\begin{theorem}\label{th:converg}
Consider the sequence $\{x^\nu\}$ generated by Algorithm \ref{algoBasic} with $\tilde f$ and $\tilde g$ such that Assumption A holds. Then, the whole sequence $\{x^\nu\}$ is contained in $K$. Furthermore,  either the sequence $\{x^\nu\}$ is unbounded or the following assertions hold:
\begin{description}
\item[\rm (i)]
at least one limit limit point $\hat x$ of $\{x^\nu\}$ is generalized stationary for problem \eqref{eq:startpro}; in particular, if  the eMFCQ holds at $\hat x$, then $\hat x$ is a KKT point for problem \eqref{eq:startpro};
\item[\rm (ii)] if, in addition, the eMFCQ holds at every limit point of $\{x^\nu\}$, under Assumption B, every limit point of  $\{x^\nu\}$ is a KKT solution for problem \eqref{eq:startpro}.
\end{description}
\end{theorem}
\proof{Proof.}
Since the starting point $x^0$ belongs to the convex set $K$, the stepsize $\gamma^\nu \leq 1$ and, by the last constraint in (P$_{x^\nu}$), $x^\nu + d(x^\nu)\in K$ for all $\nu$, it is easily seen that all points $x^\nu$ generated by the algorithm belong to $K$.
We now assume, without loss of generality, that the sequence $\{x^\nu\}$ is bounded.
Preliminarily, observe that, at each step, the solution $d(x^\nu)$ of subproblem (P$_{x^\nu}$)
is also a KKT point for (P$_{x^\nu}$). In fact, suppose that at a certain iteration $\bar \nu$, $d(x^{\bar \nu})$ does not satisfy the KKT conditions for
(P$_{x^{\bar \nu}}$). The subproblem is always feasible by construction; let us analyze the three exhaustive cases considered in Lemma \ref{th:prellem}. In case (i), Slater's condition holds for (P$_{x^{\bar \nu}}$) and $d(x^{\bar \nu})$ is a KKT point. In case (ii), $x^{\bar \nu}$ is an ES point of \eqref{eq:startpro}: hence, we would have stopped at step \ref{S.11}. In case (iii), either Slater's condition holds for (P$_{x^{\bar \nu}}$) and $d(x^{\bar \nu})$ is a KKT point, or $x^{\bar \nu}$ is a FJ point for \eqref{eq:startpro}, in which case we would have stopped at step \ref{S.11}.

 Thus, $d(x^\nu)$ is a KKT point for (P$_{x^\nu}$) and multipliers $\{\xi^\nu\}$ exist such that $\xi^\nu \in N_{\mathbb R^m_-}(\tilde g(d(x^\nu); x^\nu) - \kappa(x^\nu) e)$ and
\begin{equation}\label{eq:kktnutris}
0 \in \nabla_1 \tilde f(d(x^\nu); x^\nu) + \nabla_1 \tilde g(d(x^\nu); x^\nu) \xi^\nu + N_{\beta \mathbb B^n_\infty \cap (K - x^\nu)}(d(x^\nu)).
\end{equation}
Thanks to A1 and A4, we have
\begin{equation}\label{eq:convscobjter}
\begin{array}{rcl}
\nabla_1 \tilde f(d(x^\nu); x^\nu)\trt d(x^\nu) & = & [\nabla_1 \tilde f(d(x^\nu); x^\nu) - \nabla_1 \tilde f(0; x^\nu) + \nabla_1 \tilde f(0; x^\nu)]\trt d(x^\nu)\\[0.5em]
& \ge & c \|d(x^\nu)\|^2 + \nabla f(x^\nu)\trt d(x^\nu).
\end{array}
\end{equation}
Moreover, in view of A5, for every $i = 1, \ldots, m$,
\begin{equation}\label{eq:convscIter}
\begin{array}{rcl}
- \nabla_1 \tilde g_i(d(x^\nu); x^\nu)\trt d(x^\nu) \le \tilde g_i(0; x^\nu) - \tilde g_i(d(x^\nu); x^\nu)
\end{array}
\end{equation}
and, by A7, since $\xi^\nu$ is nonnegative, in turn,
\begin{equation}\label{eq:convscIIter}
- \xi_i^\nu \nabla_1 \tilde g_i(d(x^\nu); x^\nu)\trt d(x^\nu) \le \xi_i^\nu [\tilde g_i(0; x^\nu) - \tilde g_i(d(x^\nu); x^\nu)] = \xi_i^\nu [g_i(x^\nu) - \kappa(x^\nu)],
\end{equation}
where the equality follows observing that $\xi^\nu$ belongs to $N_{\mathbb R^m_-}(\tilde g(d(x^\nu); x^\nu) - \kappa(x^\nu) e)$.

Therefore, by \eqref{eq:kktnutris}, \eqref{eq:convscobjter} and \eqref{eq:convscIIter}, we have, for some $\zeta^\nu \in N_{\beta \mathbb B^n_\infty \cap (K - x^\nu)}(d(x^\nu)),$
\begin{equation*}\label{eq:convscIIIter}
\begin{array}{rcl}
c \|d(x^\nu)\|^2 + \nabla f(x^\nu)\trt d(x^\nu) & \le & \nabla_1 \tilde f(d(x^\nu); x^\nu)\trt d(x^\nu) = - {\xi^\nu}\trt \nabla_1 \tilde g(d(x^\nu); x^\nu)\trt d(x^\nu) - {\zeta^\nu}\trt d(x^\nu)\\[0.5em]
& \le & {\xi^\nu}\trt [g(x^\nu) - \kappa(x^\nu) e] \le {\xi^\nu}\trt [\max_i\{g_i(x^\nu)_+\} - \kappa(x^\nu)] e,
\end{array}
\end{equation*}
where the second inequality is due to $0 \in \beta \mathbb B^n_\infty \cap (K - x^\nu)$.
Therefore,  recalling definition \eqref{eq:theta},
\begin{equation}\label{eq:convscfinter}
\nabla f(x^\nu)\trt d(x^\nu) \le - c \|d(x^\nu)\|^2 +  \theta(x^\nu) \, {\xi^\nu}\trt e.
\end{equation}
We also notice that, since $d(x^\nu)$ is feasible for problem (P$_{x^\nu}$), by A5, A7 and A9,
\begin{equation}\label{eq:convscgter}
\kappa(x^\nu) \ge \tilde g_i(d(x^\nu); x^\nu) \ge \tilde g_i(0; x^\nu) + \nabla \tilde g_i(0; x^\nu)\trt d(x^\nu) = g_i(x^\nu) + \nabla g_i(x^\nu)\trt d(x^\nu).
\end{equation}
Let us now consider  the nonsmooth (ghost) penalty function already described in the introduction
\begin{equation}\label{eq:pendef}
W(x;\varepsilon) \triangleq f(x) + \frac{1}{\varepsilon} \max_i \{g_i(x)_+\},
\end{equation}
with a positive penalty parameter $\varepsilon$. This function plays a key role in the subsequent convergence analysis although it does not appear anywhere in the algorithm itself.

In the following analysis we will freely invoke some properties of function $(\bullet)_+ \triangleq \max\{0, \bullet\},$ namely $\max\{0, \alpha_1\} \le \max\{0, \alpha_2\}$ for any $\alpha_1, \, \alpha_2 \in \mathbb R$ such that $\alpha_1 \le \alpha_2,$ $\max\{0, a \, \alpha\} = a \, \max\{0, \alpha\}$ for any $\alpha \in \mathbb R$ and nonnegative scalar $a,$  and $\max\{0, \alpha_1 + \alpha_2\} \le \max\{0, \alpha_1\} + \max\{0, \alpha_2\}$ and $\max\{0, \alpha_1\} - \max\{0, \alpha_2\} \le \max\{0, \alpha_1 - \alpha_2\} $ for any $\alpha_1, \, \alpha_2 \in \mathbb R.$
We have
\begin{equation}\label{eq:nonsmmaj1ter}
\begin{array}{rcl}
W(x^{\nu + 1};\varepsilon) &-& W(x^\nu;\varepsilon)\\[5pt]
& = & f(x^\nu + \gamma^\nu d(x^\nu)) - f(x^\nu) + \frac{1}{\varepsilon} \displaystyle \Big[\max_i\{g_i(x^\nu+ \gamma^\nu d(x^\nu))_+\} \displaystyle - \max_i \{g_i(x^\nu)_+\}\Big]\\[5pt]
& \overset{(a)}{\le} & \gamma^\nu \nabla f(x^\nu)\trt d(x^\nu) + \frac{(\gamma^\nu)^2 L_{\nabla f}}{2} \|d(x^\nu)\|^2 + \frac{1}{\varepsilon} \displaystyle \Big[\max_i \{(g_i(x^\nu) + \gamma^\nu \nabla g_i(x^\nu)\trt d(x^\nu))_+\}\\[5pt]
& & - \displaystyle \max_i\{g_i(x^\nu)_+\} + \frac{(\gamma^\nu)^2  \max_i\{L_{\nabla g_i}\}}{2} \|d(x^\nu)\|^2\Big]\\[5pt]
& \overset{(b)}{\le} & \gamma^\nu \nabla f(x^\nu)\trt d(x^\nu) + \frac{1}{\varepsilon} \displaystyle \Big[\max_i \{(1 - \gamma^\nu) g_i(x^\nu)_+ + \gamma^\nu \kappa(x^\nu)\} - \displaystyle \max_i\{g_i(x^\nu)_+\}\Big]\\[5pt]
& & + \frac{(\gamma^\nu)^2}{2} (L_{\nabla f} + \frac{ \max_i\{L_{\nabla g_i}\}}{\varepsilon}) \|d(x^\nu)\|^2\\[5pt]
& \le & \gamma^\nu \nabla f(x^\nu)\trt d(x^\nu) - \frac{\gamma^\nu}{\varepsilon} \Big[\displaystyle \max_i\{g_i(x^\nu)_+\} - \kappa(x^\nu) \Big] + \frac{(\gamma^\nu)^2}{2} (L_{\nabla f} + \frac{ \max_i\{L_{\nabla g_i}\}}{\varepsilon}) \|d(x^\nu)\|^2\\[5pt]
& \le & \gamma^\nu \nabla f(x^\nu)\trt d(x^\nu) - \frac{\gamma^\nu}{\varepsilon} \, \theta(x^\nu) + \frac{(\gamma^\nu)^2}{2} (L_{\nabla f} + \frac{ \max_i\{L_{\nabla g_i}\}}{\varepsilon}) \|d(x^\nu)\|^2,
\end{array}
\end{equation}
where (a) follows applying the descent lemma to $f$ and $g_i$ for every $i = 1, \ldots, m$, with $L_{\nabla f}$ and $L_{\nabla g_i}$ being the Lipschitz moduli of $\nabla f$ and $\nabla g_i$  on the bounded set containing all iterates; (b) holds for any positive $\gamma^\nu \le 1$ since, in view of \eqref{eq:convscgter}, $\nabla g_i(x^\nu)\trt d(x^\nu) \le \kappa(x^\nu) - g_i(x^\nu)$.
Furthermore, we observe that
\begin{equation}\label{eq:nonsmmaj2ter}
\begin{array}{l}
\nabla f(x^\nu)\trt d(x^\nu) - \frac{1}{\varepsilon} \, \theta(x^\nu) \le - c \|d(x^\nu)\|^2
+ \theta(x^\nu) \, {\xi^\nu}\trt e
 - \frac{1}{\varepsilon} \, \theta(x^\nu) \le - c \|d(x^\nu)\|^2
 + (m \|\xi^\nu\|_\infty - \frac{1}{\varepsilon}) \, \theta(x^\nu),
\end{array}
\end{equation}
where the first inequality is entailed by \eqref{eq:convscfinter}.

By \eqref{eq:nonsmmaj2ter}, for any fixed $x^\nu$ and for any $\eta \in (0, 1]$, there exists $\bar \varepsilon^\nu > 0$ such that
\begin{equation}\label{eq:unifsuffdescfixter}
\nabla f(x^\nu)\trt d(x^\nu) - \frac{1}{\varepsilon} \, \theta(x^\nu) \le - \eta c \|d(x^\nu)\|^2 \qquad \forall \varepsilon \in (0, \bar \varepsilon^\nu].
\end{equation}
We now   distinguish two cases.

(I) Suppose that \eqref{eq:unifsuffdescfixter} does not hold uniformly for every $x^\nu$, that is $\eta \in (0,1]$, and a subsequence $\{x^\nu\}_{\mathcal N}$ exists, where $\mathcal{N}\subseteq\{0, 1,2, \ldots\}$ such that we can construct a corresponding subsequence $\{\varepsilon^\nu\}_{\mathcal N} \in \mathbb R_+$ with $\varepsilon^\nu \downarrow 0$ on $\mathcal N$ and
\begin{equation}\label{eq:contrter}
\nabla f(x^\nu)\trt d(x^\nu) - \frac{1}{\varepsilon^\nu} \, \theta(x^\nu) > -\eta c \|d(x^\nu)\|^2
\end{equation}
for every $\nu \in {\mathcal N}$.
For \eqref{eq:contrter} to hold, relying on \eqref{eq:nonsmmaj2ter}, the multipliers' subsequence $\{\xi^\nu\}_{\mathcal N}$ must be unbounded. Combining \eqref{eq:nonsmmaj2ter} and \eqref{eq:contrter}, we get
\begin{equation*}\label{eq:unb}
0 \le c (1 - \eta) \|d(x^\nu)\|^2 < \left(m\|\xi^\nu\|_\infty - \frac{1}{\varepsilon^\nu}\right) \, \theta(x^\nu),
\end{equation*}
and, thus, $\theta(x^\nu) > 0$ for every $\nu \in \mathcal N$.
By the previous relation and \eqref{eq:contrter}, we also have
\begin{equation}\label{eq:unb2}
\frac{1}{\varepsilon^\nu} < \frac{\nabla f(x^\nu)\trt d(x^\nu) + \eta c \|d(x^\nu)\|^2}{\, \theta(x^\nu)}.
\end{equation}
As $\varepsilon^\nu \downarrow 0$ on $\mathcal N$, the right hand side of \eqref{eq:unb2} goes to infinity: by the boundedness of the numerator,
\begin{equation}\label{eq:mustbezero}
\, \theta(x^\nu) \underset{\mathcal N}{\to} 0.
\end{equation}
Let $\hat x$ be a cluster point of the subsequence $\{x^\nu\}_{\mathcal N}$. By \eqref{eq:mustbezero}, only cases (ii) and (iii) in Lemma \ref{th:prellem} can occur at $\hat x \in K$. The existence of a $d$ as stipulated in Lemma \ref{th:prellem} (iii) would entail, by Proposition \ref{lem:boundedness}, the boundedness of the KKT multipliers $\xi^\nu$ for $\nu \in {\mathcal N}$ large enough, thus giving a contradiction. Therefore, by Lemma \ref{th:prellem} (ii), we conclude that $\hat x$ is either an ES or FJ point for \eqref{eq:startpro}.

(II) As opposed to (I), consider the case in which relation \eqref{eq:unifsuffdescfixter} holds uniformly for every $x^\nu$: that is, for any $\eta \in (0, 1]$, there exists $\bar \varepsilon > 0$ such that
\begin{equation}\label{eq:unifsuffdescter}
\nabla f(x^\nu)\trt d(x^\nu) - \frac{1}{\varepsilon} \, \theta(x^\nu) \le - \eta c \|d(x^\nu)\|^2 \enspace \forall \varepsilon \in (0, \bar \varepsilon], \;\; \forall \nu.
\end{equation}

\noindent
Combining relations \eqref{eq:nonsmmaj1ter} and \eqref{eq:unifsuffdescter}, we get
\begin{equation}\label{eq:nonsmmaj3ter}
\begin{array}{rcl}
W(x^{\nu + 1}; \tilde \varepsilon) - W(x^\nu; \tilde \varepsilon) & \le & - \gamma^\nu \eta c \|d(x^\nu)\|^2 + \frac{(\gamma^\nu)^2}{2} (L_{\nabla f} + \frac{\max_i\{L_{\nabla g_i}\}}{\tilde \varepsilon}) \|d(x^\nu)\|^2\\[5pt]
& = & -\gamma^\nu \left[\eta c -  \frac{\gamma^\nu}{2} (L_{\nabla f} + \frac{\max_i\{L_{\nabla g_i}\}}{\tilde \varepsilon})\right] \|d(x^\nu)\|^2,
\end{array}
\end{equation}
for any $\tilde \varepsilon \in (0, \bar \varepsilon]$.
Since $\lim_\nu \gamma^{\nu}=0$, there exists a positive constant $\omega$ such that, by \eqref{eq:nonsmmaj3ter}, for $\nu\ge\bar{\nu}$ sufficiently large,
\begin{equation}\label{eq:nonsmmaj2bister}
W(x^{\nu+1}; \tilde \varepsilon) -  W(x^\nu; \tilde \varepsilon) \le - \omega \gamma^\nu \|d(x^\nu)\|^2.
\end{equation}
With $W$ being bounded from below, by \eqref{eq:nonsmmaj2bister}, the sequence $\{W(x^\nu; \tilde \varepsilon)\}$ converges and
\begin{equation*}\label{eq:summable_seriester}
\lim_{\nu}\sum_{t=\bar{\nu}}^{\nu}\gamma^{t}\|d(x^t)\|^{2}<+\infty.
\end{equation*}
Therefore, since $\sum_{\nu=0}^{\infty}\gamma^{\nu}=+\infty$, we have
\begin{equation}\label{eq:liminf is 0}
\lim \inf_{\nu \to \infty} \|d(x^\nu)\| = 0.
\end{equation}
Recalling relation \eqref{eq:thetadelta}, taking the limit on a subsequence ${\mathcal N}$ such that $\|d(x^\nu)\| \underset{\mathcal N}{\to} 0$, we have $\|\nabla g(x^\nu)\|_\infty \|d(x^\nu)\| \underset{\mathcal N}{\to} 0$ and $\theta(x^\nu) \underset{\mathcal N}{\to} 0.$
Finally, let again $\hat x$ be a cluster point of subsequence $\{x^\nu\}_{\mathcal N}$. Since $\theta(x^\nu) \underset{\mathcal N}{\to} 0$ implies $\kappa (\hat x) = \max_i \{g_i(\hat x)_+\}$, cases (ii) or (iii) in Lemma \ref{th:prellem} may occur: specifically, $\hat x$ is either an ES, or a FJ, or a KKT point for \eqref{eq:startpro}. In particular, if the eMFCQ holds at $\hat x$, case (ii) in Lemma \ref{th:prellem} is ruled out and $\max_i \{g_i(\hat x)_+\}$ cannot be strictly positive; then, $\kappa (\hat x) = \max_i \{g_i(\hat x)_+\} = 0$. Furthermore, taking the limit in \eqref{eq:kktnutris}, we obtain, by A3, A4, A6-A9, KKT multipliers' boundedness and outer semicontinuity property (see Lemma \ref{th:prelresfeas} (iii)) of the normal cone mapping $N_{\beta \mathbb B^n_\infty \cap (K - \bullet)}(\bullet)$,
$$
- \nabla f(\hat x) - \nabla g(\hat x) \hat \xi \in N_{\beta \mathbb B^n_\infty \cap (K - \hat x)}(0) = \{0\} + N_{K - \hat x} (0) = N_K(\hat x),
$$
with $\hat \xi \in N_{\mathbb R^m_-}(g(\hat x) - \kappa(\hat x) e) = N_{\mathbb R^m_-}(g(\hat x))$ and where the first equality follows from Lemma \ref{th:prelresfeas} (i). In turn, $\hat x$ is a KKT point for problem \eqref{eq:startpro}. This concludes the proof of case (i).

Consider now point (ii). Note that if the eMFCQ holds at every limit point of $\{x^\nu\}$, then case (I) above cannot occur since this would contradict the last sentence before (II); hence, we are in case (II). Observe that if, instead of the weaker \eqref{eq:liminf is 0},
\begin{equation}\label{eq:lim is 0}
\lim_{\nu \to \infty} \|d(x^\nu)\| = 0
\end{equation}
holds, we can reason similarly to what done above after \eqref{eq:liminf is 0} for any convergent subsequence of $\{x^\nu\}$, and conclude that (ii) holds.
Therefore, it is enough to show that Assumption B entails \eqref{eq:lim is 0}.

Consider now the compact set containing all iterates $x^\nu$.
While $\liminf_{\nu\rightarrow\infty} \|d(x^\nu)\|=0$, suppose by contradiction that $\limsup_{\nu\rightarrow\infty} \|d(x^\nu)\|>0$. Then, there exists
$\delta>0$ such that $
 \|d(x^\nu)\|>\delta$ and $ \|d(x^\nu)\|<\delta/2
$
for infinitely many $\nu$s. Therefore, there is an infinite
subset of indices ${\cal N}$ such that, for each $\nu\in{\cal N}$,
and some $i_{\nu}>\nu$, the following relations hold:
\begin{equation}\label{eq:con1}
 \|d(x^\nu)\|<\delta/2,\hspace{6pt}\|d({x}^{i_{\nu}})\|>\delta
\end{equation}
and, if $i_\nu > \nu + 1$,
\begin{equation}
\delta/2\le\|d({x}^{j})\|\le\delta,\hspace{6pt}\nu<j<i_{\nu}.\label{eq:con2}
\end{equation}
Hence, for all $\nu\in{\cal N}$, we can write
\begin{equation}\label{eq:ineqser}
\begin{array}{rcl}
\delta/2 & < &  \|d({x}^{i_{\nu}})\|-\|d({x}^{\nu})\|
 \, \leq \, \|d({x}^{i_{\nu}}) - d({x}^{\nu})\| \overset{(a)}{\le} \theta \|{x}^{i_{\nu}}-{x}^{\nu}\|^\alpha\\[0.5em]
 & \overset{(b)}{\le} & \theta \left[\sum_{t=\nu}^{i_{\nu}-1}\gamma^{t} \|d(x^t)\|\right]^\alpha \overset{(c)}{\le} \theta \delta^\alpha\left(\sum_{t=\nu}^{i_{\nu}-1}\gamma^{t}\right)^\alpha,
\end{array}
\end{equation}
where (a) is due to Assumption B with $\alpha$ and $\theta$ positive scalars, (b) comes from the triangle inequality and the updating rule of the algorithm and in (c) we used
\eqref{eq:con2}. By \eqref{eq:ineqser} we have
\begin{equation}\label{eq:absu}
\underset{_{\nu\rightarrow\infty}}\liminf \;\; \theta \delta^\alpha\left(\sum_{t=\nu}^{i_{\nu}-1}\gamma^{t}\right)^\alpha>0.
\end{equation}
We prove next that \eqref{eq:absu} is in contradiction with the convergence of $\{W(x^\nu;\tilde \varepsilon)\}$
for any $\tilde \varepsilon \in (0, \bar \varepsilon]$, where $\bar\varepsilon$ is defined around
\eqref{eq:unifsuffdescter}. To this end, we first
show that $\|d({x}^{\nu})\|\ge\delta/4$, for sufficiently large $\nu\in{\cal N}$. Reasoning as in \eqref{eq:ineqser}, we have
\begin{equation*}
\|d({x}^{\nu+1})\| - \|d({x}^{\nu})\| \le \theta \|
{x}^{\nu+1}-{x}^{\nu}\|^\alpha \le \theta (\gamma^\nu)^\alpha \|d({x}^{\nu})\|^\alpha,
\label{eq:ineqser2}
\end{equation*}
for any given $\nu$. For  $\nu\in{\cal N}$ large enough so that $\theta ({\gamma^\nu})^\alpha(\delta/4)^\alpha<\delta/4$, suppose by contradiction that $\|
d({x}^{\nu})\|<\delta/4$; this would give $\|d({x}^{\nu+1})\|<\delta/2$ and, thus,
condition \eqref{eq:con2} (or \eqref{eq:con1}) would be violated. Then, it must be
$
\|d({x}^{\nu})\|\ge\delta/4.
$
From this, and using \eqref{eq:nonsmmaj2bister}, we have, for sufficiently large $\nu\in{\cal N}$,
\begin{equation}\label{eq:fin}
W({x}^{i_{\nu}}; \tilde \varepsilon) \le W({x}^{\nu}; \tilde \varepsilon) -\omega\sum_{t=\nu}^{i_{\nu}-1}\gamma^{t}\|d({x}^{t})\|^{2} \le W({x}^{\nu}; \tilde \varepsilon) - \omega\frac{\delta^{2}}{16}\sum_{t=\nu}^{i_{\nu}-1}\gamma^{t}.
\end{equation}
Since $\{W(x^\nu; \tilde \varepsilon)\}$ converges, as established above immediately after \eqref{eq:nonsmmaj2bister},
renumbering if necessary, relation \eqref{eq:fin} implies $\sum_{t=\nu}^{i_{\nu}-1}\gamma^{t} \to 0$, in contradiction with \eqref{eq:absu}.
This shows that \eqref{eq:lim is 0} holds and concludes the proof of the theorem. \hfill
\hfill  \Halmos \endproof

The convergence properties in Theorem \ref{th:converg} (i) are very much in the spirit of analogous results for constrained optimization where no regularity conditions are made, see for example \cite{burke1989sequential,burke1992robust,burke1989robust,Facch97}. A key difference between our approach and those in, e.g., \cite{burke1989sequential,burke1992robust,burke1989robust,Facch97} is that we do not use any penalty parameter in the algorithm. Indeed, we use
the penalty function and penalty parameter only in the proof of Theorem \ref{th:converg}, as a tool of theoretical analysis, and thus we do not need to calculate any careful penalty parameter update, allowing for convergence for the conceptually simple procedure defined above.
We believe that this ghost approach is a novelty in the literature and represents a new interesting use of penalty functions. Note that while our approach has some similarities to a classical Lyapunov function approach, it is different from it. Indeed,  while, in case (II) considered in the proof, the penalty can be viewed as a Lyapunov function for the algorithm, the analysis of  case (I) is rather different and more involved. Indeed the proof hinges on the behavior of the penalty function, of $\theta$, and of the penalty parameter and on how these quantities are connected.

\begin{remark}\label{rem: we can use it}
All the developments in the proof of Theorem \ref{th:converg} up to equation \eqref{eq:nonsmmaj3ter} are valid independent of the updating rule for the stepsize $\gamma^\nu \in (0,1]$. In the light of this observation, in the next section we invoke some of the relations in the proof of Theorem \ref{th:converg} even when stepsizes not satisfying \eqref{eq:gamma} are employed.
\end{remark}

\begin{remark}\label{rem:feasible DSM}
Algorithm \ref{algoBasic} can easily be made into a {\em feasible method}, i.e. a method that only generates
feasible iterates, if  we take $\tilde g_i$s to be UCAs, see \eqref{eq:upperb}.
As discussed in Section \ref{sec:subsec3.1}, in this case, if $x^\nu$ is feasible, then $x^{\nu+1}$ is also  feasible,
so that, if $x^0\in {\cal X}$, Algorithm \ref{algoBasic} generates only feasible iterates.
This shows that Algorithm \ref{algoBasic} contains as special cases some recent feasible
methods that were shown to be rather effective, see for example
\cite{scutari2017parallel,scutari2017parallelbis} (and also \cite{facchinei2017feasible}).

If, furthermore, we require $\tilde f$ to be an UCA for $f$, i.e.,
\begin{equation}\label{eq:uppera}
   \tilde f(d;x) \,\geq\, f(x+d),\;\;  \forall d\,\in\, K-x,
\end{equation}
we turn our scheme into a {\em Majorization-Minimization} (MM)-like method, see for example
\cite{auslender2010moving,beck2010sequential,bolte2016majorization,
hong2016unified,hunter2004tutorial,mairal2013optimization,
razaviyayn2013unified,sun2017majorization}.
In classical MM approaches the stepsize $\gamma^\nu$ is taken to be always one, while
Algorithm \ref{algoBasic} with UCAs for $f$ and $g_i$s gives a diminishing stepsize version of the method. In
Section \ref{sec:-2upper} we show that not only can we also guarantee convergence by setting the stepsize equal to one, but we can actually obtain, in this case, an iteration complexity result.
\end{remark}


\section{Iteration Complexity Analysis}\label{Sec:complexity}
We introduce some new rules for choosing the stepsize $\gamma^\nu$ at each iteration as an alternative to the diminishing one analyzed in the previous section. For these rules we are able to perform a detailed iteration complexity analysis. Our analysis is in line with recent works on this topic,  see \cite{cartis2018optimality} for an up-to-date review.
 The purpose of the iteration complexity analysis is to give a bound on the number of iterations needed by an algorithm to reach a desired level of accuracy. This bound is expressed in terms of parameters of the algorithm and some {\em problem constants}, for example Lipschitz moduli of the functions involved on a prescribed region or maximum or minimum values of the functions in the same region. This section is organized as follows. Theorem \ref{th:complbad} gives our more general complexity result for Algorithm \ref{algoCompl}; subsection \ref{subsec:meaning} explores in detail the \textcolor{black}{meaning} of the stopping criteria used in Algorithm \ref{algoCompl} \textcolor{black}{and gives the definition of $\delta$-stationary point.} The following three short subsections examine some particular scenarios in which  improved complexity bounds (or, in one case, global convergence rate) can be obtained. Finally, subsection \ref{subsec:not used} describes a variant of Algorithm \ref{algoCompl} that can be implemented and analyzed without any knowledge of any problem-related constants.

In order to perform our analysis  in this section we make the following assumptions.

\medskip \noindent {\bf{Assumption D}}
\begin{description}
\item [{\textit{D1)}}]  {\textit{the set $K$ is bounded;}}
\item [{\textit{D2)}}]   {\textit{$\nabla_1 \tilde g(\bullet; \bullet)$ is locally Lipschitz continuous on $O_d \times  O_x$.}}
\end{description}
\medskip
\noindent
Assumption D1 serves to guarantee boundedness of the iterations and is made for simplicity of presentation.
In Section \ref{sec:bounded}  we shall discuss some alternative assumptions that make the iterates belong to a compact set defined by means of possibly known quantities, as required in order to perform a complexity analysis, see in particular Table \ref{tab2} and the surrounding comments. As the discussion pertains to the algorithms presented in both the previous and current sections, we found that a detailed analysis of this condition is best deferred in order to not complicate the formal presentation of the results, as the insight involved is essentially modular, separate from the main ideas of analysis here and in Section~\ref{Sec:convergence}.

\noindent Assumption D2, instead, depends essentially on the choice of $\tilde g$ and therefore is not an assumption on the problem itself, but a condition on our algorithmic choices. Clearly, since $g$ has a locally Lipschitz gradient, D2 is always satisfied if we take as $\tilde g$ the linearization of $g$.

From now on, we employ some problem dependent constants: we collect their definitions in Table \ref{tab1} for the reader's convenience.

We observe that if the eMFCQ holds everywhere in $K$, all generalized stationary solutions are KKT points for problem \eqref{eq:startpro} and, as in  classical SQP methods, the norm of the direction $d(x^\nu)$ is a natural stationarity measure, see Theorem \ref{th:stopping criteria}. But if the eMFCQ is not valid at every point in $K$, we cannot rely solely on $\|d(x^\nu)\|$ to monitor progress towards stationarity, since the problem may admit KKT points but also FJ and ES solutions. For this reason, we use in combination $\|d(x^\nu)\|$ and $\theta(x^\nu)$ as measures of stationarity.
We observe that, actually, $\|d(x^\nu)\|$ and $\theta(x^\nu)$ are linked to each other in view of the following relation, which is due to \eqref{eq:thetadelta}:
\begin{equation}\label{eq:thetadelta2}
\theta(x^\nu) \le \|\nabla g(x^\nu)\|_\infty \|d(x^\nu)\| \le L \|d(x^\nu)\|,
\end{equation}
where $L \triangleq \max_{(d, x)} \{\|\nabla_1 \tilde g(d;x)\|_\infty \, | \, (d, x) \in \beta \mathbb B_\infty^n \times K\}$. However, there is no reverse implication and thus the two functions $\|d(x^\nu)\|$ and $\theta(x^\nu)$ must be suitably combined to provide reliable stopping criteria.
The  effect of monitoring both $\|d(x^\nu)\|$ and $\theta(x^\nu)$ on the outcome of the algorithm is analyzed
in detail in Section \ref{subsec:meaning}.

\begin{table}
\TABLE
{Problem dependent constants \label{tab1}}
{\begin{tabular}{l@{\hskip 2.5em}l}
\hline
\\
$\beta \in \mathbb R_+$ & user-set constant in the definition of \eqref{eq:p_k} \\[0.5em]
$\eta \in (0,1]$ & user-set constant\\[0.5em]
$\lambda \in (0,1)$ & user-set constant in the definition of $\kappa$, see \eqref{eq:cap} \\[0.5em]
$c \in \mathbb R_+$ & modulus of strong convexity of $\tilde f(\bullet;x)$, see Assumption A1 \\[0.5em]
$B$ & $\max_{x \in K} \|\nabla f(x)\| \beta + \eta c \beta^2$\\[0.5em]
$f^m$ & $\min_{x\in K} f(x)$\\[0.5em]
$g^M_+$ & $\max_{x\in K} \max_i\{ g_i(x)_+\}$\\[0.5em]
$L$ & $\max_{(d,x) \in \beta \mathbb B^n_\infty \times K} \|\nabla_1 \tilde g(d; x)\|_\infty$\\[0.5em]
$L_{\nabla f}$ & Lipschitz modulus of $\nabla f$ on $K$\\[0.5em]
$L_{\nabla g_i}$ & Lipschitz modulus of $\nabla g_i$ on $K$\\[0.5em]
$L_{\nabla \tilde f}$ & Lipschitz modulus of $\nabla_1 \tilde f(\bullet; \bullet)$  on $\beta {\mathbb B}_\infty^n\times K$\\[0.5em]
$L_{\nabla \tilde g}$ & Lipschitz modulus of $\nabla_1 \tilde g(\bullet; \bullet)$ on $\beta {\mathbb B}_\infty^n\times K$ \\[0.5em]
\end{tabular}}
{}
\end{table}

To derive complexity results, we consider first Algorithm \ref{algoCompl} with a piecewise constant choice of stepsizes. By this we mean that Algorithm \ref{algoCompl} starts with a certain $\gamma^{-1}$ and keeps it fixed until a certain test is met; when this happens, the stepsize is reduced to a new prescribed value and then kept fixed until possibly the test is met again, and so on.
 We underline that the only difference between this scheme and Algorithm \ref{algoBasic} is in the rules for choosing $\gamma^\nu$ at each iteration and, of course, in the presence of suitable stopping criteria: specifically, the steps \ref{S.1} and \ref{S.7} correspond to the previous Algorithm \ref{algoBasic}, while everything in between, from \ref{S.2} to \ref{S.6}, is aimed at deciding whether to decrease the stepsize $\gamma^\nu$ and  whether we should terminate (note that Algorithm \ref{algoBasic}, which is aimed at an asymptotic analysis, does not contain any practical stopping criterion).

\medskip

{\setlength{\skiptext}{10pt}
\setlength{\skiprule}{5pt}

\IncMargin{2em}
\begin{algorithm}[H]
\KwData{\small{$\delta > 0$, $\eta \in (0,1]$,  $x^{0} \in K$, $T^{-1} \in \left(0,\frac{2 \max_i \{L_{\nabla g_i}\}}{\max \{L_{\nabla f}, \eta c\}}\right]$, $\gamma^{-1} =  \frac{T^{-1} \eta c}{2 \max_i\{L_{\nabla g_i}\}}$, $\nu \longleftarrow 0$\;}}
\Repeat{
\nlset{(S.1)}compute $\kappa(x^\nu)$, the solution $d(x^\nu)$ of problem (P$_{x^\nu}$) and $\theta(x^\nu)$\; \label{S.1} \smallskip

{\nlset{(S.2)} \If{$\|d(x^{\nu})\| \le \delta$}{{\bf{stop}} and {\bf{return}} $x_\delta = x^\nu$\;} \label{S.2}
{\nlset{(S.3)} \uIf{\small{$\nabla f(x^\nu)\trt d(x^\nu) + \eta c \|d(x^\nu)\|^2 > 0$} {\bf and} \small{$T^{\nu-1} > \frac{\theta(x^\nu)}{\nabla f(x^\nu)\trt d(x^\nu) + \eta c \|d(x^\nu)\|^2}$} \label{S.3}}{\nlset{(S.4)} \uIf{$\theta(x^{\nu})  \le \delta$\label{S.4}}{{\bf{stop}} and {\bf{return}} $x_\delta = x^\nu$\;}\Else{\nlset{(S.5)} set \small{$\gamma^\nu = \frac{T^\nu \eta c}{2 \max_i\{L_{\nabla g_i}\}}$, where $T^\nu = \frac{1}{2} \frac{\theta(x^\nu)}{\nabla f(x^\nu)\trt d(x^\nu) + \eta c \|d(x^\nu)\|^2}$\;}\label{S.5}}}
\Else{\nlset{(S.6)} set {\small{$T^\nu = T^{\nu - 1}$}} and {\small{$\gamma^\nu = \gamma^{\nu - 1}$}}\;}\label{S.6}}}
\nlset{(S.7)} set $x^{\nu+1}=x^{\nu}+\gamma^{\nu}d(x^\nu)$, $\nu\longleftarrow\nu+1$\;\label{S.7}}{}
\caption{\label{algoCompl}Modified Algorithm for \eqref{eq:startpro}}
\end{algorithm}}
\DecMargin{2em}

\medskip

\noindent
{
We first note that the value of $T^{-1}$ guarantees that $\gamma^{-1} \leq 1$. Also, the variable  $T^\nu$ is introduced just for notational purposes, in order to make the statement of the algorithm and the proof of Theorem \ref{th:complbad} easier to follow. The tests  we must perform to decide whether to reduce the stepsize are very simple and involve quantities that are readily available once the direction finding subproblem (P$_{x^\nu}$) has been solved. The following theorem provides the announced complexity result in this general case. For simplicity of presentation we assume  $\delta \leq 1$. This is by no means necessary but avoids the necessity of  complicating the statement by considering uninteresting cases.

\begin{theorem}\label{th:complbad}
Let $\{x^\nu\}$ be the sequence generated by Algorithm \ref{algoCompl} under Assumptions A, C1 and D and suppose that $\delta \leq 1$.
Then,  in at most $\mathcal{O}(\delta^{-4})$
iterations, Algorithm \ref{algoCompl} stops either at step \emph{\ref{S.2}} or at step \emph{\ref{S.4}}; more precisely, the maximum number of iterations is given by the maximum between the expressions \eqref{eq:rare case} and \eqref{eq:numberiterations1}.
\end{theorem}

\proof{Proof.}
Suppose that Algorithm \ref{algoCompl} performs $N$ iterations without stopping\footnote{We consider an iteration completed when we reach \ref{S.7}.}.
We first count how many times $\gamma^\nu$ can be updated in step \ref{S.5} of the algorithm: let
\[{\mathcal{I}} \triangleq \{0<\nu_i\leq N \, | \, T^{\nu_i} \, \text{and} \, \gamma^{\nu_i} \, \text{are updated in \ref{S.5}}\}\, \cup\,  \{0\}
\]
be the set of iterations' indices $\nu$ (in increasing order) at which the need to modify $\gamma^\nu$ and $
T^\nu$ emerges, union iteration 0. Therefore, for example, if we update $T$ and $\gamma$ in \ref{S.5} at iterations
3, 4 and 8, we have ${\mathcal{I}} = \{ \nu_0 = 0, \nu_1 =3, \nu_2 = 4, \nu_3 =8\}$; note that we always have by definition $\nu_0 =0$ and that the set $\cal I$ does not include repeated indices.
We show  that $\mathcal{I}$ has finite cardinality.
If $\nu_i\neq 0$ belongs to ${\mathcal{I}}$, we have
\begin{equation}\label{eq:upd}
T^{\nu_i} = \frac{1}{2} \frac{\theta(x^{\nu_i})}{\nabla f(x^{\nu_i})\trt d(x^{\nu_i}) + \eta c \|d(x^{\nu_i})\|^2},
\end{equation}
and the procedure did not stop  at step \ref{S.4}: thus, $\theta(x^{\nu_i}) >  \delta$ and \eqref{eq:upd} entails
\begin{equation}\label{eq:upd2}
T^{\nu_i} > \frac{ \delta}{2 B},
\end{equation}
 with $B \triangleq \max_x \left\{\|\nabla f(x)\| \beta + \eta c \beta^2 \, | \, x \in K\right\} \ge \nabla f(x^\nu)\trt d(x^\nu) + \eta c \|d(x^\nu)\|^2$.
By the updating rule in \ref{S.5}, we also have $T^{\nu_i} \le \frac{T^{-1}}{2^{i}}$; thus, in view of \eqref{eq:upd2}, $\frac{ \delta}{2 B} < T^{\nu_i} \le \frac{T^{-1}}{2^{i}}$, so that
\begin{equation*}\label{eq:outitcount}
i < \log_2 \frac{T^{-1} 2 B}{ \delta}.
\end{equation*}
Therefore, if we do not stop, i.e. if  $\theta(x^\nu)>\delta$ for all iterations up to $N - 1$, the cardinality of ${\mathcal{I}}$, i.e., the times  $\gamma^\nu$ is reduced, is  at most  $\ceil*{\log_2 \frac{T^{-1} 2 B}{ \delta}}$.

Let us set  $I\triangleq |{\cal I}| -1$; with this convention note that the largest element in $\cal I$ is $\nu_I$.
Counting from $\nu_i\in {\cal I}\setminus\{$last element in ${\cal I} \}$, let now $N_i$
be the number of  iterations in which
 $\gamma^\nu$ remains unchanged:
 $T^\nu = T^{\nu_i}$ and $\gamma^\nu = \gamma^{\nu_i}$ for every $\nu \in  \{\nu_i, \ldots, \nu_i + N_i\}$. In other words, $N_i$ is the number of iterations after $\nu_i$ in which step \ref{S.5} is not reached; in the example given above where ${\mathcal{I}} = \{ \nu_0 = 0, \nu_1 =3, \nu_2 = 4, \nu_3 =8\}$, we have
 $N_0 = 2$,  $N_1 = 0$,  $N_2 = 3$. Therefore $\nu_i +N_i$ is simply the last iteration after $\nu_i$ before $\gamma$ and $T$ are updated. The last index $N_{I}$ is defined, with the same rationale, as the number of iterations performed after $\nu_{I}$, before we reach the iteration where we stop. Considering the example  above, and supposing that we stop at iteration $11$, we have $N_3 = 2$.

We observe that, by virtue of the condition in step \ref{S.3} and the updating rule in step \ref{S.5} or \ref{S.6}, $T^\nu$ is non increasing. Hence, again by the updating rule in \ref{S.5} or \ref{S.6}, since $\gamma^{-1} = \frac{T^{-1} \eta c}{2 \max_i\{L_{\nabla g_i}\}}$, also $\gamma^\nu$ is non increasing.
Moreover,by the definitions of  $T^{-1}$ and $\gamma^{-1}$,  on the one hand, $\eta c - \frac{\gamma^\nu}{2} L_{\nabla f} \ge \eta c - \frac{\gamma^{-1}}{2} L_{\nabla f} \ge \eta c - \frac{\eta c}{2}$, while, on the other hand, $-\frac{\gamma^\nu}{2} \frac{\max_i \{L_{\nabla g_i}\}}{T^\nu} = - \frac{\eta c}{4}$. Thanks to the previous relations, we have for every $\nu$
\begin{equation}\label{eq:conststepsmall}
\eta c - \frac{\gamma^\nu}{2} \left(L_{\nabla f} + \frac{\max_i \{L_{\nabla g_i}\}}{T^\nu}\right) \ge \frac{\eta c}{4},
\end{equation}
and, in turn, by \eqref{eq:nonsmmaj1ter},
\begin{equation}\label{eq:preldisc}
W(x^{\nu+1}; T^\nu) - W(x^\nu; T^\nu) \le \gamma^\nu \left[\nabla f(x^\nu)\trt d(x^\nu) - \frac{\theta(x^\nu)}{T^{\nu}} + \frac{3\eta c}{4} \|d(x^\nu)\|^2\right],
\end{equation}
where we took $\varepsilon^\nu = T^\nu$. We now distinguish two cases. If the condition in step \ref{S.3} is satisfied and $\gamma$ is updated in \ref{S.5},
\begin{equation}\label{eq:preldisc2}
\nabla f(x^\nu)\trt d(x^\nu) - \frac{\theta(x^\nu)}{T^{\nu}} = - \nabla f(x^\nu)\trt d(x^\nu) - 2 \eta c \|d(x^\nu)\|^2 \le - \eta c \|d(x^\nu)\|^2,
\end{equation}
where the inequality follows from the first condition in \ref{S.3}. If, on the contrary, $\gamma$ need not be reduced,
\[
\nabla f(x^\nu)\trt d(x^\nu) + \eta c \|d(x^\nu)\|^2 \le 0 \quad {\text{or}} \quad T^{\nu} \le \frac{\theta(x^\nu)}{\nabla f(x^\nu)\trt d(x^\nu) + \eta c \|d(x^\nu)\|^2},
\]
and, again, relation \eqref{eq:preldisc2} is easily seen to hold. Therefore, in view of \eqref{eq:preldisc} and \eqref{eq:preldisc2}, we get for every $\nu$,
\begin{equation}\label{eq:desccomplbis}
W(x^{\nu + 1}; T^{\nu}) - W(x^\nu; T^{\nu}) \le - \gamma^{\nu} \, \frac{\eta c}{4} \, \|d(x^\nu)\|^2.
\end{equation}
Note that $N=\sum_{i \in {\cal I}}(N_i+1)$, since the algorithm did not stop until iteration $N$,  $\|d(x^\nu)\| > \delta$ for all iterates up to $N - 1$.
Therefore,  recalling definition \eqref{eq:pendef} with $\varepsilon^\nu = T^\nu$, and observing that for every $\nu \in \{\nu_i, \ldots, \nu_i + N_i\}$,  $\nu_i \in {\cal I}$,   $\gamma^\nu$ is not reduced and $T^\nu = T^{\nu_i}$, we get
\begin{equation}\label{eq:majbad}
\begin{array}{rcl}
\displaystyle \delta^2 N & = & \displaystyle \sum_{i =0}^{I} \delta^2 (N_i + 1) < \displaystyle
\sum_{i=0 }^I \sum_{\nu = \nu_i}^{\nu_i + N_i}\|d(x^\nu)\|^2 \le \sum_{i =0}^I \displaystyle \frac{W(x^{ \nu_i}; T^{\nu_i}) - W(x^{\nu_i + N_i + 1}; T^{\nu_i})}{\gamma^{\nu_i} \, \frac{\eta c}{4}}\\[5pt]
& \le & \frac{1}{\gamma^{\nu_I} \frac{\eta c}{4}}\Big[f(x^{0}) - f(x^{\nu_I + N_I + 1}) + \frac{1}{T^{0}} \max_i \{g_i(x^{0})_+\}\\[5pt]
& &- \frac{1}{T^{\nu_I + N_I + 1}} \max_i \{g_i(x^{\nu_I + N_I + 1})_+\} + \displaystyle \sum_{i=1}^I \left(\frac{1}{T^{\nu_i}} - \frac{1}{T^{\nu_{i - 1}}} \right) \max_j \{g_j(x^{\nu_i})_+\}  \Big],
\end{array}
\end{equation}
where the second inequality is due to \eqref{eq:desccomplbis} while, observing that $\gamma^{\nu_I} \le \gamma^{\nu_i}$, the last inequality is valid as a result of a telescopic series argument since $\nu_i + N_i + 1 = \nu_{i + 1}$. Note that it is understood that if $I=0$ the last summation in \eqref{eq:majbad} has no terms.
Letting $g_+^M \triangleq \max_x \{\max_i \{g_i(x)_+\} \, | \, x \in K\}$ and $f^m \triangleq \min_x \{f(x) \, | \, x \in K\}$, we distinguish two cases: (i) step \ref{S.5} has never been reached, i.e. $T$ has never been diminished; (ii) case (i) did not occur. In case (i), observing that $I =0$, by \eqref{eq:majbad}, the algorithm stops after at most
\begin{equation}\label{eq:rare case}
\ceil*{\frac{8}{(\eta c)^2 T^{-1}} \max_i\{L_{\nabla g_i}\} [f(x^0) - f^m + \frac{1}{T^{-1}} \max_i \{g_i(x^0)_+\}]\frac{1}{\delta^2}}
\end{equation}
iterations. In case (ii), by \eqref{eq:majbad} we can write instead
\begin{equation}\label{eq:badmajor}
\begin{array}{rcl}
\displaystyle \delta^2 N &  \overset{(a)}{<} & \frac{1}{\gamma^{\nu_I} \frac{\eta c}{4}}\left(f(x^0) - f^m + \frac{1}{T^{0}} g_+^M - \frac{1}{T^{0}} g_+^M + \frac{1}{T^{\nu_I}} g_+^M \right) \overset{(b)}{=} \frac{8}{(\eta c)^2 T^{\nu_I}} \max_i\{L_{\nabla g_i}\} (f(x^0) - f^m + \frac{1}{T^{\nu_I}} g_+^M)
\end{array}
\end{equation}
where (a), since $T^{\nu_i} \le T^{\nu_{i - 1}}$, follows again from the summation of a telescopic series,  (b) is due to the updating rule for $\gamma^\nu$ in \ref{S.5} at iteration $\nu_I$. In turn, taking into account that since we updated $T$ at least once, we have
\[
T^{\nu_I} = \frac{1}{2} \frac{\theta(x^{\nu_I})}{\nabla f(x^{\nu_I})\trt d(x^{\nu_I}) + \eta c \|d(x^{\nu_I})\|^2}  > \frac{ \delta}{2 B},
\]
and, in turn,
\[
\delta^2 N < \frac{16B}{(\eta c)^2 \delta} \max_i\{L_{\nabla g_i}\} (f(x^0) - f^m + \frac{2B}{\delta} g_+^M),
\]
thus meaning that the procedure halts in at most
\begin{equation}\label{eq:numberiterations1}
\ceil*{\frac{16B}{(\eta c)^2} \max_i\{L_{\nabla g_i}\} \left[\frac{f(x^0) - f^m}{\delta^3} + \frac{2 B \, g_+^M}{ \delta^4}\right]}
\end{equation}
iterations. If $\delta \leq 1$, this gives an overall complexity of ${\mathcal{O}}(\delta^{-4})$.
\hfill  \Halmos
\endproof


\subsection{On the meaning of the stopping criteria at \ref{S.2} and \ref{S.4}}\label{subsec:meaning}
The following theorem elucidates the meaning of the stopping criteria in steps \ref{S.2} and \ref{S.4}. This result and the ensuing discussion show that \ref{S.2} and \ref{S.4} guarantee that the algorithm stops in a finite number of iterations once a $\delta-$stationary point has been reached. To simplify the proof, we assume that $\delta < \min\{1, \beta\}$; this is very sensible since on the one hand we are mainly interested in what happens when $\delta$ is ``small" and, on the other hand, $\beta$ is chosen by the user and is intended to be ``large", $\beta$ being simply a safeguard on the maximum length of the direction $d(x^\nu)$.

  Preliminarily, we recall that  the KKT conditions at a point  $x^\nu \in K$ for problem \eqref{eq:startpro} can be rewritten as
\begin{equation}\label{eq:projkkt}
\begin{array}{rcl}
\left\Vert P_K\left(x^\nu - \frac{\nabla f(x^\nu) + \nabla g(x^\nu) \xi^\nu}{1 + \|\xi^\nu\|}\right) - x^\nu\right\Vert = 0, \enspace \max_i \left|g_i(x^\nu) \frac{\xi_i^\nu}{1 + \|\xi^\nu\|}\right| = 0, \enspace \max_i \{g_i(x^\nu)_+\} = 0,\\[5pt]
\end{array}
\end{equation}
where $P_K$ denotes the projection on the closed convex set $K$ and $\xi^\nu \ge 0$ are suitable multipliers. We also recall that $\theta$ is a stationarity measure for the violation-of-the-constraint problem \eqref{eq:feaspro}: $\theta(x^\nu) = 0$ if and only if $x^\nu$ is stationary for \eqref{eq:feaspro}, see (ii) Proposition \ref{th:theta}.
\begin{theorem}\label{th:stopping criteria}
	Let Assumptions A, C1 and D hold, and consider $\delta < \min\{1, \beta\}$. If Algorithm \ref{algoCompl}
	\begin{itemize}
		\item[(i)]  stops at step \ref{S.2}, $x^\nu$ is either infeasible almost stationary for the violation-of-the-constraints problem, i.e.,
\begin{equation}\label{eq:deg}
\max_i \{g_i(x^\nu)_+\} > \frac{L}{\lambda} \delta, \enspace 0 < \theta(x^\nu) \le L \, \delta,
\end{equation}
or it is a scaled-KKT point, i.e.,
\begin{equation}\label{eq:scaledKKT}
\begin{array}{c}
\max_i \{g_i(x^\nu)_+\} \le \frac{L}{\lambda} \, \delta\\[5pt]
\left\Vert P_K\left(x^\nu - \left[\frac{1}{1 + \|\xi^\nu\|} \nabla f(x^\nu) + \nabla g(x^\nu) \frac{\xi^\nu}{1 + \|\xi^\nu\|}\right]\right) - x^\nu\right\Vert \le (2+L_{\nabla \tilde f}+L_{\nabla \tilde g}) \, \delta,\\[10pt] 
\max_i \left|g_i(x^\nu) \, \frac{\xi_i^\nu}{1 + \|\xi^\nu\|}\right| \le \frac{1+\lambda}{\lambda} \, L \, \delta, 
\end{array}
\end{equation}
for some $\xi^\nu \ge 0$, or either an ES or a FJ point;

  \item[(ii)]  stops at step \ref{S.4}, $x^\nu$ is either an ES or a FJ point, or it is infeasible almost stationary for the violation-of-the-constraints problem, i.e.,
\begin{equation}\label{eq:deg2}
\max_i \{g_i(x^\nu)_+\} > 0, \enspace 0 < \theta(x^\nu) \le \delta.
\end{equation}

	\end{itemize}	
\end{theorem}
Before proving the theorem, some comments are in order.
Theorem \ref{th:stopping criteria} shows that Algorithm \ref{algoCompl} stops with a KKT, FJ or ES solution, or a point that, at least, satisfies either \eqref{eq:deg}, or \eqref{eq:scaledKKT} or \eqref{eq:deg2}. This outcome is in line with many recent results in the literature. Relations \eqref{eq:scaledKKT}, and \eqref{eq:deg} and \eqref{eq:deg2}
are similar
to classical conditions such as (i) and (ii) in \cite[Theorem 2.9]{cartis2017corrigendum}, (3.27) and (3.26)
(respectively) in \cite[Theorem 3.8]{cartis2019evaluation}, or (10) and (9) (respectively) in \cite[Theorem 4.5]
{cartis2018optimality}. Specifically, we obtain the scaled-type conditions \eqref{eq:scaledKKT}: we refer the interested reader to \cite[Section 2.1]{cartis2017corrigendum}, but also
\cite{birgin2016evaluation,cartis2013evaluation,cartis2019evaluation} for rather exhaustive discussions on this
point. On the other hand, the ``degenerate'' cases \eqref{eq:deg} and \eqref{eq:deg2} indicate, although in slightly different ways, that a stationarity
condition for the violation-of-the-constraint problem is approximately satisfied at an infeasible point.
In order to get more
insight into the meaning of  the stopping criteria, we discuss, in the same spirit as the analysis
in \cite{birgin2016evaluation}, what happens when $\delta$ goes to zero with fixed initial data. Thus,
suppose we have a sequence $\{x_k\}$ each point of which satisfies at least one of \eqref{eq:deg},
\eqref{eq:scaledKKT} or
\eqref{eq:deg2} for a sequence of values $\delta^k \downarrow 0$: in fact, we recall that, in view of Theorem \ref{th:complbad}, for every $k$ the algorithm stops, providing $x_k$, either at step \ref{S.2} or at step \ref{S.4} in a finite number of iterations $N = N^k$ (which is obviously nondecreasing with respect to $k$). Accordingly, let $I = I^k$ be the corresponding number of times $T^\nu$ and $\gamma^\nu$ have been reduced, apart from iteration $0$. Moreover, since $\{x_k\}$ is contained in $K$, it is bounded and therefore we can assume, without loss of generality, that it converges to a point $\bar x\in K$.

Suppose first that $x_k$ satisfies the scaled-KKT condition \eqref{eq:scaledKKT} for every $k \in {\mathcal K} \subseteq \{0, 1,2, \ldots \}$ for  some $\cal K$. Passing to the limit in \eqref{eq:scaledKKT}, if the corresponding sequence $\xi_k$ is bounded, $\bar x$ is a KKT point of problem \eqref{eq:startpro}. If, instead, $\xi_k$ is unbounded, $\bar x$ must be a FJ point, since it is feasible by the first inequality in \eqref{eq:scaledKKT} and the eMFCQ cannot hold there, otherwise the sequence $\xi_k$ would be bounded by Proposition \ref{th:keyemfcq} (i) and Proposition \ref{lem:boundedness}.

Suppose now that one of the two degenerate cases
 \eqref{eq:deg} or \eqref{eq:deg2} occurs at each $x_k$ with $k \in {\mathcal K}$.
 We show that, for both cases, $\bar x$ is a point where the eMFCQ does not hold and therefore it is either a FJ or an ES point.  Let the algorithm stop at step \ref{S.2} providing $x_k$ that satisfies \eqref{eq:deg} for all $k \in {\mathcal K}$, and assume by contradiction that the eMFCQ holds at $\bar x$. It follows that  $d(x_{k}) \to d(\bar x)= 0$ and $\theta(x_{k}) \to \theta(\bar x) = 0$, due to the condition $d(x_{k}) \le \delta^k$ for every $k \in {\mathcal K}$ and to the continuity relative to $K$ of function $d(\bullet)$ (on a neighborhood of $\bar x$, see Proposition \ref{th:keyemfcq}) and $\theta(\bullet)$ (see Proposition \ref{th:theta}), respectively. Besides, relying on Lemma \ref{th:prellem} (iii), $\bar d \in \rho \mathbb B_\infty^n \cap (K - \bar x)$ exists such that $\tilde g(\bar d;\bar x) < 0$. Thus, by  the continuity relative to $K$ of the set-valued mapping $K - \bullet$ at $\bar x$ (see Lemma \ref{th:prelresfeas} (ii)), there exists $d_k \in \rho \mathbb B_\infty^n \cap (K - x_k)$ such that, for every $k \in \mathcal K$ sufficiently large, $\tilde g(d_k; x_{k}) < 0$. In turn, $\min_d \left\{\max_i\{\tilde g_i(d; x_{k})_+\} \, | \, \|d\|_\infty \le \rho, \, d \in K - x_{k}\right\} = 0$, $\kappa(x_{k}) = (1 - \lambda) \max_i \{g_i(x_{k})_+\}$, and $\theta(x_{k}) = \lambda \max_i \{g_i(x_{k})_+\} \le L \delta^k$ in contradiction with $\max_i \{g_i(x_{k})_+\} > \frac{L}{\lambda} \delta^k$. Therefore the eMFCQ does not hold at $\bar x$.

Suppose now that the algorithm stops at step \ref{S.4} for all $k \in {\mathcal K}$ and assume by contradiction that the eMFCQ holds at $\bar x$.  We have, without loss of generality, $N^{k} > N^{k-1}$ for every $k \in {\mathcal K}$ and $\theta(x_{k}) \to\theta(\bar x) = 0$, due to the condition $\theta(x_{k}) \le \delta^k$ for every $k \in {\mathcal K}$ and to the continuity relative to $K$ of function $\theta(\bullet)$. Furthermore, it holds $I^{k} \ge I^{k-1} + 1$ for every $k \in {\mathcal K}$ and, in turn, $T^{N^k} \downarrow 0$ on ${\mathcal K}$, since $T^{N^k} = T^{\nu_{I^k}} \le \frac{T^{-1}}{2^{I^{k}}}$ for every $k \in {\mathcal K}$. If  the eMFCQ holds at $\bar x$, for any $k \in {\mathcal K}$ sufficiently large, $d(x_{k})$ is a KKT point for (P$_{x_{k}}$) by Proposition \ref{th:cqlocal} and, in turn, by
\eqref{eq:nonsmmaj2ter}, we get
\begin{equation}\label{eq:contrcomplff}
\begin{array}{rcl}
\nabla f(x_{k})\trt d(x_{k}) - \frac{1}{T^{N^k}} \theta(x_{k}) + \eta c \|d(x_{k})\|^2 \le
\left(m \|\xi^{N^k}\|_\infty - \frac{1}{T^{N^k}} \right) \theta(x_{k}).
\end{array}
\end{equation}
Thanks to the local boundedness of the set of KKT multipliers and because $T^{N^k} \downarrow 0$ on $\mathcal K$, eventually the right hand side of \eqref{eq:contrcomplff} is nonpositive, in contradiction to the condition $\nabla f(x_{k})\trt d(x_{k}) + \eta c \|d(x_{k})\|^2 > 0$ and $T^{N^k} > \frac{\theta(x_{k})}{\nabla f(x_{k})\trt d(x_{k}) + \eta c \|d(x_{k})\|^2}$ for every $k \in \mathcal K$ in \ref{S.3}.  Therefore the eMFCQ does not hold at $\bar x$.
\textcolor{black}{All the above discussion motivates us to define a point at which Algorithm \ref{algoCompl} stops a $\delta-$(generalized) stationary point.
\begin{definition}\label{def:delta approx}
A point  generated by the algorithm is a $\delta-$(generalized) stationary point if it is either a scaled-KKT point satisfying  \eqref{eq:scaledKKT} or an infeasible approximate stationary point for the   violation-of-the-constraints-problem satisfying
\eqref{eq:deg} or
\eqref{eq:deg2}.
\end{definition}
}

It may also be interesting to remark that if the eMFCQ holds at every point in $K$, \eqref{eq:deg} and \eqref{eq:deg2} cannot occur if $\delta$ is small enough (see the discussion above), and $\xi^\nu$ that, we shall see in the proof below, are the multipliers of the direction finding subproblems (P$_{x^\nu}$), are bounded by Propositions \ref{th:keyemfcq} (i) and \ref{lem:boundedness}. In turn, this means that the algorithm stops at \ref{S.2} with a point $x^\nu$ approximately satisfying the KKT conditions for \eqref{eq:startpro} with $\xi^\nu$ being nothing else but approximate KKT multipliers (see \cite[Section 2.1]{cartis2017corrigendum} for further details).

\medskip

\proof{Proof of Theorem \ref{th:stopping criteria}.}
(i) Suppose first that the algorithm stops because $\|d(x^\nu)\| \le \delta$. Regardless of the validity of the constraint qualification, $d(x^\nu)$, which certainly satisfies the Fritz-John conditions, may satisfy or not the KKT conditions for the subproblem (P$_{x^\nu}$). We now distinguish two cases, remarking that the following results hold whatever the choice of $\gamma^\nu$.

	(I) If $d(x^\nu)$ does not satisfy the KKT conditions for subproblem (P$_{x^\nu}$), in view of Proposition \ref{th:keyemfcq}, $x^\nu$ does not satisfy the eMFCQ and, thus, is either an ES or a FJ point.
	
	(II) If, on the contrary, $d(x^\nu)$ satisfies the KKT conditions for subproblem (P$_{x^\nu}$), letting $\xi^\nu \in N_{\mathbb R_-^m}(\tilde g(d(x^\nu); x^\nu) - \kappa(x^\nu) e)$, we get the following relation which is equivalent to \eqref{eq:kktnutris} that still holds with $N_{\beta \mathbb B^n_\infty} = \{0\}$ because $\|d(x^\nu)\| \le \delta$:
\begin{equation}\label{eq:projopt}
	x^\nu + d(x^\nu) = P_K\left(x^\nu + d(x^\nu) - \frac{\nabla_1 \tilde f(d(x^\nu); x^\nu) + \nabla_1 \tilde g(d(x^\nu); x^\nu) \xi^\nu}{1 + \|\xi^\nu\|}\right).
\end{equation}
Let us bound now the terms in relations \eqref{eq:projkkt}. As for the gradient of the Lagrangian-related condition for problem \eqref{eq:startpro}, we have the following bound:
\begin{equation*}\label{eq:meas1}
\begin{array}{rcl}
\left\Vert P_K\left(x^\nu - \frac{\nabla f(x^\nu) + \nabla g(x^\nu) \xi^\nu}{1 + \|\xi^\nu\|}\right) - x^\nu\right\Vert & = &  \left\Vert P_K\left(x^\nu - \frac{\nabla f(x^\nu) + \nabla g(x^\nu) \xi^\nu}{1 + \|\xi^\nu\|}\right) - [x^\nu + d(x^\nu)] + d(x^\nu)\right\Vert\\[5pt]
& \overset{(a)}{=} & \left\Vert d(x^\nu) + P_K\left(x^\nu - \frac{\nabla f(x^\nu) + \nabla g(x^\nu) \xi^\nu}{1 + \|\xi^\nu\|}\right)\right.\\[5pt]
& & \left. - P_K\left(x^\nu + d(x^\nu) - \frac{\nabla_1 \tilde f(d(x^\nu); x^\nu) + \nabla_1 \tilde g(d(x^\nu); x^\nu) \xi^\nu}{1 + \|\xi^\nu\|}\right)\right\Vert\\[5pt]
& \overset{(b)}{\le} & \|d(x^\nu)\| + \left\Vert-d(x^\nu) + \frac{\nabla_1 \tilde f(d(x^\nu); x^\nu) - \nabla_1 \tilde f(0; x^\nu)}{1 + \|\xi^\nu\|}\right.\\[5pt]
& & \left. + \frac{\nabla_1 \tilde g(d(x^\nu); x^\nu) \xi^\nu - \nabla_1\tilde g(0; x^\nu) \xi^\nu}{1 + \|\xi^\nu\|}\right\Vert\\[5pt]
& \overset{(c)}{\le} & (2 + L_{\nabla  \tilde f} + L_{\nabla  \tilde g}) \, \|d(x^\nu)\|
\end{array}
\end{equation*}
where (a) follows from \eqref{eq:projopt}, (b) holds thanks to A4, A9 and since the projection mapping is nonexpansive, and (c) is due to C1 and D2. As for the complementarity conditions, consider $\bar \imath \in \{1, \ldots, m\}$ such that $|g_{\bar \imath}(x^\nu) \xi_{\bar \imath}^\nu| = \max_i |g_i(x^\nu) \xi_i^\nu|$ with $\tilde g_{\bar \imath}(d(x^\nu); x^\nu) = \kappa(x^\nu)$, otherwise $\xi_{\bar \imath}^\nu = 0$. Note that, if $g_{\bar \imath}(x^\nu) \ge 0$, $g_{\bar \imath}(x^\nu) \le \max_i \{g_i(x^\nu)_+\}$,
whereas if $g_{\bar \imath}(x^\nu) < 0$,
\[
|g_{\bar \imath}(x^\nu)| = -g_{\bar \imath}(x^\nu) - \tilde g_{\bar \imath}(d(x^\nu); x^\nu) + \tilde g_{\bar \imath}(d(x^\nu); x^\nu) \le \nabla_1 \tilde g_{\bar \imath}(d(x^\nu); x^\nu)\trt d(x^\nu), 	
\]
where the inequality is due to \eqref{eq:convscIter} and $-\tilde g_{\bar \imath}(d(x^\nu); x^\nu) = -\kappa(x^\nu) \le 0$. Overall,
\begin{equation*}\label{eq:compconbis}
\max_i |g_i(x^\nu) \, \xi_i^\nu| \le (\max_i \{g_i(x^\nu)_+\} + L \|d(x^\nu)\|) \|\xi^\nu\|. 
\end{equation*}
In turn, if $\max_i \{g_i(x^\nu)_+\} \le \frac{L}{\lambda}\delta$, then
\begin{equation*}\label{eq:compcondter}
\max_i \left|g_i(x^\nu) \, \frac{\xi^\nu}{1 + \|\xi^\nu\|}\right| \le \frac{1 + \lambda}{\lambda} \, L \, \delta.	
\end{equation*}
If, on the contrary, $\max_i \{g_i(x^\nu)_+\} > \frac{L}{\lambda} \delta$, nonetheless, by \eqref{eq:thetadelta2} we have $\theta(x^\nu) \le L \, \delta$.

(ii) In order to exit at step \ref{S.4}, either $x^\nu$ is an ES or a FJ point, or $\theta(x^\nu)$ must be strictly positive. In fact, under the eMFCQ, by \eqref{eq:convscfinter}, if $\theta(x^\nu) = 0$, then $\nabla f(x^\nu)\trt d(x^\nu) \le - \eta c \|d(x^\nu)\|^2$ and the first condition in step \ref{S.3} does not hold. In turn, for $\theta(x^\nu)$ to be strictly positive, we must have $\max_i \{g_i(x^\nu)\} > 0$.
\hfill  \Halmos
\endproof

\subsection{$O(\delta^{-3})$ complexity with constant stepsize if a feasible starting point is known}\label{sec:-3}
If a feasible starting point is available, then by choosing  a sufficiently small initial $T^{-1}$ or, correspondingly, a sufficiently small initial stepsize $\gamma^{-1}$, the iteration complexity of Algorithm \ref{algoCompl} can be reduced to ${\mathcal{O}}(\delta^{-3})$.  Actually, it turns out that in this case, as well as in the cases analyzed in the next two subsections, the stepsize is never reduced, so that the updating step of Algorithm \ref{algoCompl} actually becomes a {\em fixed stepsize iteration}
\begin{equation}\label{eq:fixed stepsize}
x^{\nu + 1} = x^\nu + \bar \gamma d(x^\nu).
\end{equation}
A reduction of the iteration complexity when a feasible point is available  seems rather sensible because if we start with a feasible point, we have already solved the feasibility problem which is a part of the constrained optimization. Nevertheless, it was in principle not clear that our algorithm could take advantage of this fact, since the search for feasibility and that for optimality are combined in a single step, unlike typical methods
designed for strong complexity results for constrained nonconvex problems that use two distinct phases.

\begin{corollary} Assume the same setting of Theorem \ref{th:complbad}, fix a prescribed tolerance $\delta$ and set, according to this value,  $T^{-1} = \min\{
\frac{\delta}{B}, \frac{2 \max_i \{L_{\nabla g_i}\}}{\max\{L_{\nabla f}, \eta c\}}\}$. If the starting
point $x^0$ is feasible, then, in at most
\[
\ceil*{\frac{8}{(\eta c)^2} \max_i\{L_{\nabla g_i}\} \max \left\{B, \; \frac{\max\{L_{\nabla f}, \, \eta c\}}{2\max_i\{L_{\nabla g_i}\}}\right\} (f(x^0) - f^m) \frac{1}{\delta^3}}
\]
 iterations, Algorithm \ref{algoCompl} stops
 either at step \emph{\ref{S.2}} or at step \emph{\ref{S.4}}. Furthermore, the stepsize is never updated and is constant throughout the algorithm.
\end{corollary}

\proof{Proof.} We use  the same notation and terminology
introduced in the proof of Theorem \ref{th:complbad}.
We first observe that Algorithm \ref{algoCompl} never updates $\gamma^\nu$ and $T^\nu$. Indeed, suppose that the test in \ref{S.3} is met for the first time at iteration $\nu$.
The claim follows noting that if the condition in  \ref{S.3} is verified, then
\[
\frac{ \delta}{B} \geq T^{-1} =  T^{\nu-1} > \frac{\theta(x^\nu)}{\nabla f(x^\nu)\trt d(x^\nu) + \eta c \|d(x^\nu)^2\|} \ge \frac{\theta(x^\nu)}{B},
\]
so that $\theta(x^\nu) \le \delta$ and the algorithm stops. Hence, step \ref{S.5} is never reached and the stepsize is never updated. Looking back at the corresponding case (i) in Theorem \ref{th:complbad}, in view of \eqref{eq:rare case} and recalling that $\delta \leq 1$, the procedure is shown to stop after the claimed number of iterations, at worst.
\hfill\Halmos
\endproof
Note the somewhat unusual feature that  algorithmic choices, i.e. $T^{-1}$, are linked  to the desired accuracy.

\subsection{$O(\delta^{-2})$  complexity with constant stepsize if a feasible starting point is known and  upper approximations are used}\label{sec:-2upper}

Suppose again that a feasible starting point is available and, in addition, assume that UCAs for the $g_i$s are
used, (see Section \ref{sec:subsec3.1} and \eqref{eq:upperb} in particular). Then, not only can
we get $O(\delta^{-2})$ complexity, but, differently from the previous subsection, there is no dependence of
$T^{-1}$ on $\delta$. Also in this case it turns out that the stepsize need not be updated, and
the algorithm reduces to the fixed stepsize scheme \eqref{eq:fixed stepsize}. Furthermore, if we are in an MM
setting, i.e. if we choose an UCA also for $f$ (see Remark \ref{rem:feasible DSM} and \eqref{eq:uppera}), we can
take the fixed stepsize to be one, provided some minimal assumptions on $\tilde f$ are satisfied.

\begin{corollary} Assume the same setting of Theorem \ref{th:complbad}. If the starting point $x^0$ is feasible
and the $\tilde g_i$s are upper convex approximations for the $g_i$s, then, in at most
\[
\ceil*{\frac{8}{(\eta c)^2 T^{-1}} \max_i\{L_{\nabla g_i}\} [f(x^0) - f^m]\frac{1}{\delta^2}}
\]
iterations, Algorithm \ref{algoCompl} stops either at step \emph{\ref{S.2}} or at step \emph{\ref{S.4}}. Furthermore, the stepsize is never updated, is constant throughout the algorithm progress, and can be set equal to one provided that $\eta c \geq L_{\nabla f}$.
\end{corollary}

\proof{Proof.} The algorithm only produces feasible iterates, see the discussion after \eqref{eq:upperb}. Therefore, we have $\theta(x^\nu) =0$ for all $\nu$. As a consequence, step \ref{S.5} is never reached and the stepsize is never updated: in fact, if the test in \ref{S.3} is met, then the algorithm immediately stops at \ref{S.4} since $\theta(x^\nu) = 0$. Then, reasoning again as in case (i) in Theorem \ref{th:complbad}, thanks to \eqref{eq:rare case}, we see that the algorithm stops after at most the claimed number of iterations.
Suppose further that $\eta c \geq L_{\nabla f}$, then it is easy to see from the instructions in {\bf  Data} that we
can choose $\gamma^{-1}=1$.\hfill \Halmos
\endproof

We remark that it is easy to show that the condition $\eta c \geq L_{\nabla f}$ implies that $\tilde f$ is an UCA
and therefore the requirement in the corollary imposes that we use not any arbitrary UCA, but
only  UCAs that additionally
satisfy $\eta c \geq L_{\nabla f}$.
At the same time, in standard MM algorithms it is usually possible  to  show convergence  with a unitary stepsize
without requiring $\eta c \geq L_{\nabla f}$,  or similar assumptions.
 But we must observe that
the constants $\eta$ and $c$ are algorithmic choices and therefore the condition   $\eta c \geq L_{\nabla f}$ can
always be enforced.
For example, if analogously to what done in \eqref{eq:UCAg}, we set
\begin{equation}\label{eq:UCAf}
 \tilde f(d; x) = f(x) + \nabla f(x)\trt d + \frac{c}{2}\| d\|^2,
 \end{equation}
it is enough to choose $c$ so that  $\eta c \geq L_{\nabla f}$.
Additionally, and more importantly,
 the condition $\eta c \geq L_{\nabla f}$ is needed here in order to establish for the first time, as far as we are
 aware of, the iteration complexity for an MM method.
 Our iteration complexity complements the convergence rate obtained in \cite{bolte2016majorization}. In that paper, assuming a Kurdyka-\L{}ojasiewicz property plus other technical conditions, the authors show that, under suitable constraint qualifications, that we do not require, the whole sequence produced by an MM method converges to a KKT point $x^\infty$ and give expressions for the convergence rate of $\|x^\nu - x^\infty\|$.

\subsection{A $O(\delta^{-2})$ global convergence rate when eMFCQ holds} \label{sec:-2}

If the eMFCQ holds at every point in $K$, then we can prove that Algorithm \ref{algoCompl} has a {\em global convergence rate} of $O(\delta^{-2})$. Once again, under suitable assumptions, one can show that in Algorithm \ref{algoCompl} the stepsize is never updated, so that the algorithm reduces to the fixed stepsize iteration \eqref{eq:fixed stepsize}.

\begin{corollary}\label{cor:fixed} Assume the same setting of Theorem \ref{th:complbad} and, in addition, suppose that the eMFCQ holds at every point in $K$.
If we choose $T^{-1}$ and, correspondingly, $\gamma^{-1}$ sufficiently small (as will be specified in the proof, see also the comments below), being $M$ an upper bound on the norm of multipliers for the subproblems (P$_{x^\nu}$), in at most \begin{equation}\label{eq:rare bis}
\ceil*{\frac{8 m M}{(\eta c)^2 } \max_i\{L_{\nabla g_i}\} [f(x^0) - f^m + m M \max_i \{g_i(x^0)_+\}]\frac{1}{\delta^2}}
\end{equation} iterations, Algorithm \ref{algoCompl} stops at step \emph{\ref{S.2}}. Furthermore, the stepsize is never updated and is constant throughout the algorithm.
\end{corollary}
 Note that since we never reach \ref{S.4}, the only stopping criterion actually used is the one based on $\|d(x^\nu)\|$ in \ref{S.2}, in accordance with what happens in classical SQP-type methods when constraint qualifications are assumed to hold everywhere.

\proof{Proof of Corollary \ref{cor:fixed}.}
We first recall that, thanks to the eMFCQ, by Propositions \ref{th:keyemfcq} (i) and \ref{lem:boundedness}, taking into account the compactness of $K$, the norm of multipliers $\xi^\nu$ of the subproblems (P$_{x^\nu}$) is bounded from above by some constant $M$.
 By \eqref{eq:nonsmmaj2ter}, which, in view of the eMFCQ, is still valid because it is derived from the optimality conditions for subproblem (P$_{x^\nu}$), we have
\begin{equation}\label{eq:contrcompl}
\begin{array}{rcl}
\nabla f(x^\nu)\trt d(x^\nu) - \frac{1}{T^\nu} \theta(x^\nu) + \eta c \|d(x^\nu)\|^2 \le \left(m \|\xi^\nu\|_\infty - \frac{1}{T^\nu} \right) \theta(x^\nu) \le \left(m M - \frac{1}{T^{-1}} \right) \theta(x^\nu) \le 0
\end{array}
\end{equation}
if $T^{-1} \le 1/mM$, and, in turn, for all $\nu$ it never happens that $\nabla f(x^\nu)\trt d(x^\nu) + \eta c \|d(x^\nu)\|^2 > 0$ and $T^\nu > \frac{\theta(x^\nu)}{\nabla f(x^\nu)\trt d(x^\nu) + \eta c \|d(x^\nu)\|^2}$ in \ref{S.3}, and therefore, \ref{S.4} and \ref{S.5} are never reached. Looking back at the proof of Theorem \ref{th:complbad}, we then see that only case (i) therein can occur and, setting, for example, $T^{-1}=\frac{1}{mM}$ in \eqref{eq:rare case}, the algorithm stops at worst after the claimed number of iterations.
 \hfill \Halmos
\endproof
The fixed stepsize iteration \eqref{eq:fixed stepsize} is valid provided that $\bar \gamma \le \frac{ \eta c}{2mM \max_i\{L_{\nabla g_i}\}}$. Finally, we remark that the bound given by \eqref{eq:rare bis} is different  from those seen so far, in that it depends on the usually unknown quantity $M$. The bound given by \eqref{eq:rare bis} should therefore be regarded as a {\em global convergence rate} (see the introduction).

\subsection{Problem constants are not used}\label{subsec:not used}
The implementation of Algorithm \ref{algoCompl} requires the use of some of the problem constants in Table \ref{tab1}. Hence, the question arises whether we can modify the algorithm to avoid the use of potentially difficult to compute constants, while retaining complexity results similar to those in Theorem \ref{th:complbad}. The answer is positive, at the price of a ``small amount'' of additional function evaluations. Moreover, differently from all previous developments, we must make a numerical, although simple,  use of the penalty function $W$. Observe that in  Algorithm \ref{algoCompl} the problem constants are used to set some initial values in {\bf Data} and, more critically, in \ref{S.5}. Referring to the proof of  Theorem \ref{th:complbad}, the updating of $\gamma^\nu$ in \ref{S.5} guarantees condition \eqref{eq:desccomplbis}, i.e. the sufficient decrease of the (ghost) penalty function. But, at a more basic level, this sufficient decrease condition can always be reached if the step $\gamma^\nu$ is sufficiently small. So, one could choose at each iteration the stepsize $\gamma^\nu$ so as to guarantee that the sufficient decrease condition \eqref{eq:desccomplbis} is satisfied. This can be accomplished without any knowledge of the problem constants; we only need to know the user-set quantities $c$ and $\eta$ as shown in Algorithm \ref{algoComplb}.

\medskip

{\setlength{\skiptext}{10pt}
\setlength{\skiprule}{5pt}

\IncMargin{2em}
\begin{algorithm}
\KwData{\small{$\delta > 0$, $\eta \in (0,1]$, $x^{0} \in K$, $T^{-1} >0$, $\gamma^{-1}  =1$, $\nu \longleftarrow 0$\;}}
\Repeat{
\nlset{(S.1)}compute $\kappa(x^\nu)$, the solution $d(x^\nu)$ of problem (P$_{x^\nu}$) and $\theta(x^\nu)$\; \label{S.1b} \smallskip

{\nlset{(S.2)} \If{$\|d(x^{\nu})\| \le \delta$}{{\bf{stop}} and {\bf{return}} $x_\delta = x^\nu$\;} \label{S.2b}

{\nlset{(S.3)} \uIf{\small{$\nabla f(x^\nu)\trt d(x^\nu) + \eta c \|d(x^\nu)\|^2 > 0$} {\bf and} \small{$T^{\nu-1} > \frac{\theta(x^\nu)}{\nabla f(x^\nu)\trt d(x^\nu) + \eta c \|d(x^\nu)\|^2}$} \label{S.3b}}{\nlset{(S.4)} \uIf{$\theta(x^{\nu})  \le \delta$\label{S.4b}}{{\bf{stop}} and {\bf{return}} $x_\delta = x^\nu$\;}\Else{\nlset{(S.5)} set \small{$T^\nu = \frac{1}{2} \frac{\theta(x^\nu)}{\nabla f(x^\nu)\trt d(x^\nu) + \eta c \|d(x^\nu)\|^2}$\;}\label{S.5b}}}
\Else{\nlset{(S.6)} set {\small{$T^\nu = T^{\nu - 1}$}}\;}\label{S.6b}}}

{\nlset{(S.7)} \While{$W(x^\nu+\gamma^\nu d(x^\nu);T^\nu)-W(x^\nu;T^\nu) > -\gamma^\nu \frac{\eta c}{4}\|d(x^\nu)\|^2$} {set $\gamma^\nu \longleftarrow  \frac12 \gamma^\nu$ \;}\label{S.7b}}
\nlset{(S.8)} set $x^{\nu+1}=x^{\nu}+\gamma^{\nu}d(x^\nu)$, $\nu\longleftarrow\nu+1$\;\label{S.8b}}{}

\caption{\label{algoComplb}Algorithm for \eqref{eq:startpro} without constants}
\end{algorithm}}
\DecMargin{2em}

\medskip

\noindent
In {\bf Data} we no longer need to set the initial $T$ and $\gamma$ to some small values that depend on problem constants. Indeed, whatever the initial values, it is the algorithm itself that sets them to the appropriate quantities. In  Algorithm \ref{algoCompl}, updating the stepsize at \ref{S.5} makes \eqref{eq:desccomplbis} satisfied at all subsequent iterations, until the {\bf if} section at \ref{S.3} is possibly re-entered. In Algorithm \ref{algoComplb} instead, we do not have such a guarantee, and thus we perform the ``line-search'' in \ref{S.7} at each iteration. The following theorem shows that Algorithm \ref{algoComplb} needs an amount of iterations which is similar (likely smaller, see the comments after the proof) to that required by Algorithm \ref{algoCompl}. However, while for Algorithm \ref{algoComplb} this quantity is also equal to the number of function and constraints evaluations, we now may need some extra objective and constraint function evaluations, as detailed next.

\begin{theorem}\label{th:complbadb}
Let $\{x^\nu\}$ be the sequence generated by Algorithm \ref{algoComplb} under Assumptions A, C1 and D and suppose that $\delta \leq 1$. Then, in at most $\mathcal{O}(\delta^{-4})$
iterations, Algorithm \ref{algoComplb} stops either at step \emph{\ref{S.2b}} or at step \emph{\ref{S.4b}}; more precisely, the maximum number of iterations is given by the maximum between the expressions \eqref{max number iterations ff} and \eqref{eq:second complexity}.
Moreover, the algorithm needs a number of objective and constraint function evaluations that is equal
to the number of iterations plus at most ${\mathcal{O}}\left(\log_2(\delta^{-1})\right)$ further evaluations, with the precise expression of this additional number of evaluations  given by \eqref{eq:extra evaluations}.
\end{theorem}

\proof{Proof.} The proof is a variant of that of Theorem \ref{th:complbad}, to which we refer for notation and terminology. Suppose that Algorithm \ref{algoComplb} performs $N$ iterations without stopping.
We first count how many times $T^\nu$ can be updated in step \ref{S.5b}: let
\[{\mathcal{I}} \triangleq \{0<\nu_i\leq N \, | \, T^{\nu_i}  \; \text{is updated in \ref{S.5b}}\}\, \cup\,  \{0\}
\]
be the set of iterations' indices $\nu$ (in increasing order) at which we need to modify $T^\nu$, union iteration 0.
Repeating {\textit{verbatim}} the first part in the proof of Theorem \ref{th:complbad}, one can show that $\mathcal{I}$ has finite cardinality and, if $\nu_i \in {\cal I}$ then
\begin{equation*}\label{eq:outitcountb}
i < \log_2 \frac{T^{-1} 2 B}{ \delta}.
\end{equation*}
Define now  $I$ and $N_i$ as  in the proof of Theorem \ref{th:complbad}. Clearly $T^\nu = T^{\nu_i}$ for every $\nu \in  \{\nu_i, \ldots, \nu_i + N_i\}$. Following the same line of reasoning as in the proof of Theorem \ref{th:complbad}, one can readily show that
\begin{equation}\label{eq:condbound}
\nabla f(x^\nu)\trt d(x^\nu) - \frac{\theta(x^\nu)}{T^{\nu}} \le - \eta c \|d(x^\nu)\|^2,
\end{equation}
for every $\nu$, which, in turn, by \eqref{eq:nonsmmaj1ter}, implies
\begin{equation}\label{eq:desccomplb}
W(x^{\nu + 1}; T^{\nu}) - W(x^\nu; T^{\nu}) \le - \gamma^{\nu} \left[\eta c-\frac{\gamma^\nu}{2}(L_{\nabla f}
+\frac{\max_i\{L_{\nabla g_i}\}}{T^{\nu}})\right] \|d(x^\nu)\|^2,
\end{equation}
 where we took $\varepsilon^\nu = T^\nu$. We now note that for every $\nu \in \{\nu_i, \ldots, \nu_i + N_i\}$, $\nu_i \in \cal I$, we have $\gamma^\nu \geq G(\delta)$, with
\begin{equation}\label{eq:def G}
G(\delta) \triangleq \min \left\{\frac{1}{2}, \; \frac{3 \eta c}{4} \frac{\delta}{L_{\nabla f} \delta + 2B \max_i \{L_{\nabla g_i}\}}\right\}.
\end{equation}
Indeed, this is trivial for $\nu_0=0$, since we assumed $\gamma^{-1}=1$ and $G(\delta)\leq \frac{1}{2}$. Suppose by contradiction that $0\neq \nu_i \in \cal I$ and $\gamma^\nu < G(\delta)$. By the definition \eqref{eq:def G} of $G(\delta)$, if we set $\gamma^\nu < 2G(\delta)$, we get, recalling \eqref{eq:desccomplb} and $T^{\nu_i} > \frac{\delta}{2B}$, $W(x^{\nu + 1}; T^{\nu_i}) - W(x^\nu; T^{\nu_i}) \le - \gamma^\nu \frac{\eta c}{4} \|d(x^\nu)\|^2,$ i.e. the test at \ref{S.7b} is surely not satisfied if $\gamma^\nu < 2G(\delta)$. This, in turn, contradicts $\gamma^\nu < G(\delta)$, since it shows that in the loop \ref{S.7} we should have stopped at the previous iterate of the cycle.
 Therefore, taking into account that $\gamma^\nu$ is obtained at \ref{S.7b} after a certain number (possibly zero) of halvings  of the current value of the stepsize, at each iteration $\gamma^\nu \geq G(\delta)$ for every $\delta$. We conclude that $\gamma^\nu$, globally,  needs to be halved no more than $\log_2 \, \frac{\gamma^{-1}}{G(\delta)}$ times: hence, after at most
 \begin{equation}\label{eq:extra evaluations}
 \log_2 \left(\gamma^{-1} \max \left\{2, \; \frac{4}{3 \eta c}\left(L_{\nabla f} \delta + 2 B \max_i\{L_{\nabla g_i}\}\right)\right\} \frac{1}{\delta}\right)
 \end{equation}
halvings, one achieves the sought decrease condition
\begin{equation}\label{eq:desccomplq}
W(x^{\nu + 1}; T^{\nu}) - W(x^\nu; T^{\nu}) \le - \gamma^{\nu} \frac{\eta c}{4} \|d(x^\nu)\|^2 \le - G(\delta) \frac{\eta c}{4} \|d(x^\nu)\|^2,
\end{equation}
for every $\nu$. Recalling definition \eqref{eq:pendef}, and similarly to \eqref{eq:majbad},
\begin{equation}\label{eq:majbadb}
\begin{array}{rcl}
\displaystyle \delta^2 N & = & \displaystyle \sum_{i = 0}^I \delta^2 (N_i + 1) < \displaystyle \sum_{i=0}^I \sum_{\nu = \nu_i}^{\nu_i + N_i}\|d(x^\nu)\|^2 \le \sum_{i = 0}^I \displaystyle \frac{W(x^{ \nu_i}; T^{\nu_i}) - W(x^{\nu_i + N_i + 1}; T^{\nu_i})}{G(\delta) \, \frac{\eta c}{4}}\\[5pt]
& \le & \frac{1}{G(\delta) \frac{\eta c}{4}}\Big[f(x^{\nu_0}) - f(x^{\nu_I + N_I + 1}) + \frac{1}{T^{\nu_0}} \max_i \{g_i(x^{\nu_0})_+\}\\[5pt]
& &- \frac{1}{T^{\nu_I + N_I + 1}} \max_i \{g_i(x^{\nu_I + N_I + 1})_+\} + \displaystyle \sum_{i=1}^I \left(\frac{1}{T^{\nu_i}} - \frac{1}{T^{\nu_{i - 1}}} \right) \max_j \{g_j(x^{\nu_i})_+\}  \Big],
\end{array}
\end{equation}
where the second inequality is due to \eqref{eq:desccomplq}, while the last inequality is valid as a result of a telescopic series argument since $\nu_i + N_i + 1 = \nu_{i + 1}$. Setting $g_+^M \triangleq \max_x \{\max_i \{g_i(x)_+\} \, | \, x \in K\}$ and $f^m \triangleq \min_x \{f(x) \, | \, x \in K\}$, as in the proof of Theorem \ref{th:complbad}, we distinguish two cases: (i) step \ref{S.5b} has never been reached, thus $I=0$ and, by \eqref{eq:majbadb}, the algorithms stops after at most
\begin{equation}\label{max number iterations ff}
\ceil*{\frac{4}{\eta c}\max \left\{2, \;\frac{4}{3 \eta c} \left(L_{\nabla f}\delta + 2B \max_i \{L_{\nabla g_i}\}\right)\right\}  \left(f(x^0) - f^m + \frac{1}{T^{-1}} \max_i \{g_i(x^0)_+\} \right)\frac{1}{\delta^3}}
\end{equation}
iterations, in view of the definition of $G$ and $\delta \le 1$. If case (i) did not occur, by \eqref{eq:majbadb} we can write
\begin{equation}\label{eq:badmajorbis}
\displaystyle \delta^2 N   <  \frac{1}{G(\delta) \frac{\eta c}{4}}\left(f(x^0) - f^m + \frac{1}
{T^{\nu_0}} g_+^M - \frac{1}{T^{\nu_0}} g_+^M + \frac{1}{T^{\nu_I}} g_+^M \right),
\end{equation}
where the inequality follows, recalling that  $T^{\nu_i} \le T^{\nu_{i - 1}}$,  from the summation of a telescopic series.
In turn, since
\[
T^{\nu_I} = \frac{1}{2} \frac{\theta(x^{\nu_I})}{\nabla f(x^{\nu_I})\trt d(x^{\nu_I}) + \eta c \|d(x^{\nu_I})\|^2}  > \frac{ \delta}{2 B},
\]
by \eqref{eq:badmajorbis}, the procedure halts in at most
\begin{equation}\label{eq:second complexity}
\ceil*{\frac{4}{\eta c} \max \left\{2, \;\frac{4}{3 \eta c } \left(L_{\nabla f}\delta + 2B \max_i \{L_{\nabla g_i}\}\right)\right\} \left[\frac{f(x^0) - f^m}{\delta^3} + \frac{2 B \, g_+^M}{\delta^4}\right]}
\end{equation}
iterations. Since $\delta \le 1$, one can take the overall complexity to be of ${\mathcal{O}}(\delta^{-4})$. \hfill  \Halmos
\endproof
It is interesting to compare the worst-case bounds \eqref{eq:numberiterations1} (for Algorithm \ref{algoCompl}) and \eqref{eq:second complexity} (for Algorithm \ref{algoComplb}). It is clear that, at least for a small $\delta$ (see the definition of $G$), the bound \eqref{eq:second complexity} is approximatively $\frac{2}{3}$ of the bound \eqref{eq:numberiterations1}. This better behaviour of Algorithm \ref{algoComplb} has a simple explanation. The steps used in Algorithm \ref{algoComplb} are generally larger than those used in Algorithm \ref{algoCompl}, where problem constants are used to define a ``pessimistic'' step-length. In Algorithm \ref{algoComplb}, instead, local information is gathered through the line-search in \ref{S.7b} that permits the definition of a stepsize better adapted to the problem. Algorithm \ref{algoComplb} also has the additional merit of not requiring the knowledge of the problem constants. We pay a price for this better result in that the algorithm is marginally more complex, requires a numerical use of the penalty function and calls for additional  objective function and constraint evaluations that may increase the computational effort. However, note that this increase is negligible when $\delta$ is small, since the additional number of function evaluations is ${\mathcal{O}}\left(\log_2(\delta^{-1})\right)$, implying that the overall order of function evaluations is maintained to be ${\mathcal{O}}\left(\delta^{-4}\right)$.
\begin{remark}
It is easy to see that the results in Sections \ref{sec:-3}-\ref{sec:-2} for Algorithm \ref{algoCompl} can be extended to Algorithm \ref{algoComplb}, but we do not pursue this for lack of space.
\end{remark}

\begin{remark}\label{rem:compl} While the analysis in this section is about SQP-type approaches, a few other complexity bounds are available in the literature for different methods using first-order information in the context of nonconvex, constrained optimization (see the introduction). A direct comparison of all these results is difficult, since different assumptions and, in some cases, different concepts of (approximate) generalized stationary points are called for; furthermore, the various methods differ markedly in the overall structure (penalty vs. Phase I - Phase II vs. SQP schemes) and in the algorithmic computational effort required at each iteration, not to mention that in some cases results are given for equality constraints only. With this in mind, here we try to briefly highlight the main features of
\cite{cartis2011evaluation,cartis2013evaluation,birgin2016evaluation,cartis2017corrigendum}. The analysis provided in \cite[Section 3.2]{cartis2011evaluation} is for a  penalty-based approach and gives the ``first worst-case global evaluation bounds for constrained optimization when both the objective and the constraints are allowed to be nonconvex". The main similarity with our analysis is in the  use of a nondifferentiable penalty function, which however is employed according to a classical double-loop scheme: at each outer iteration a penalty parameter is chosen and then the penalty function is minimized inexactly using a trust-region-like method for which complexity results are provided. This should be contrasted with our ``ghost" use of penalties where the penalty function and the penalty parameter are not used in the algorithm itself. Another difference is in the subproblems to be solved at each iteration.
Our method follows the standard SQP approach and the subproblems should be regarded as  (simple) approximations of the original problem and as such do not involve any penalty parameter. The subproblems in \cite{cartis2011evaluation}, instead, aim at approximating the penalty function and, as such, necessarily include the penalty parameter.  On the one hand, avoiding the use of the penalty parameter in the subproblems is a favourable numerical feature, we believe, on the other hand approximating directly the penalty function, as done in \cite{cartis2011evaluation}, nicely avoids the issue of the feasibility of subproblems, since the minimization of the penalty function is unconstrained.
Putting together a judicious analysis of the parameter-updating scheme and of the inner penalty minimization, the authors of \cite{cartis2011evaluation} can then give estimates for the number of iterations necessary to reach an approximate generalized stationary point. If, during the minimization process, the penalty parameter grows unbounded, a complexity of ${\mathcal{O}}(\delta^{-5})$ is obtained. If, on the other hand, an upper bound, which in principle is unknown in advance, for the penalty parameter exists, a condition that should be interpreted as a constraint qualification, then, using the terminology of the present paper, a {\em convergence rate} of ${\mathcal{O}}(\delta^{-2})$ is obtained. Note that  (a) the issue of the boundedness of the iterates is not dealt with (for our algorithm see also next section), and (b) the objective function $f$ is assumed to be bounded from below on $\mathbb R^n$.

The Phase I - Phase II approach in
\cite[Sections 4 and 5]{cartis2013evaluation} is for equality constrained problems, and relies on target-following and (inner) trust-region techniques.
Both phases are based on a cubic regularization method involving the use of second order derivatives, with subproblems to be solved at each iterations that are potentially expensive. However, a remarkable complexity bound of ${\mathcal{O}}(\delta^{-3/2})$, matching the one for the cubic regularization method in the unconstrained case, is achieved.

The high-level scheme put forward in \cite{birgin2016evaluation} for nonconvex problems with (equality and) inequality constraints falls also within a Phase I - Phase II (of target-following-type) framework: the unconstrained nonlinear minimization problems  to be solved in each phase are  assumed to be dealt with by some minimization algorithm with known complexity guarantees. Once this algorithm is given, the analysis derives bounds that range between ${\mathcal{O}}(\delta^{-3})$ and ${\mathcal{O}}(\delta^{-5})$ according to the choice of an algorithmic parameter, with the better complexity corresponding to ``weaker" notions of almost generalized stationary points. It has to be remarked that, as a major departure from all other works, in \cite{birgin2016evaluation} an emphasis is given to {\em unscaled} KKT conditions, i.e. to an approximate notion of stationarity that does not depend on the magnitude of the multipliers involved.

Finally, in \cite{cartis2017corrigendum}, a Phase I - Phase II (of target-following-type) method is presented,
which resorts again to an inner trust-region approach in both phases. Assuming the objective function to be upper and lower bounded on the feasible set, and the gradient of the objective and the Jacobian of the constraint functions to be Lipschitz continuous on $\mathbb R^n$ and on a suitable extended neighbourhood of the feasible region, respectively, the algorithm is proven to reach an approximate generalized stationary point in at most ${\mathcal{O}}(\delta^{-2})$ iterations. These results
are currently the most advanced for a first order Phase I - Phase II method and, remarkably, the bounds obtained there match the best result for first-order unconstrained minimization methods.
\end{remark}

\section{Boundedness of iterates}\label{sec:bounded}

Boundedness of the sequence generated by an SQP-type method is a difficult issue. With a few earlier exceptions,
see e.g. \cite{Facch97}, this topic probably came  to a wider attention only with the important paper
\cite{solodov2009global}, that motivated researchers to look better into this issue, see
\cite{auslender2010moving,auslender2013extended,bolte2016majorization,liu2011sequential} (with the latter
reference dealing only with equality constraints, however). In our framework, generating
an unbounded sequence is a natural possibility that cannot and should not  be excluded in principle, since
we do not make  any standard  assumption such as feasibility, existence of an optimal solution, or regularity of
the constraints; quoting from \cite{bolte2016majorization}, where a similar possibility is considered, ``The divergence property... is a positive result, a
convergence result, which does not correspond to a failure of the method but rather to the absence of minimizers
in a given zone."  To clarify this point consider
\[
\begin{array}{cl}
\underset{x}{\mbox{minimize}} & x^2 \\
\mbox{s.t.} & e^x \le 0,
\end{array}
\]
which is an infeasible convex problem and has no ES, FJ or KKT solutions. Nevertheless, we can apply
one of the algorithms studied in this paper to it and the only sensible outcome is ``an attempt
to minimize infeasibility" with the generation of an unbounded sequence. Indeed, if the sequence generated by
the algorithm were bounded, every limit point should be critical, but since there are no critical points, the
sequence must necessarily be unbounded. And yet, also in the spirit of the works mentioned at the beginning of
this subsection, it is of course of great interest to see under what conditions we can guarantee the boundedness
of the sequence $\{ x^\nu\}$ for the algorithms presented in this paper. Below we analyze this
point identifying some settings where boundedness can be guaranteed.
This discussion,
which in no way tries to be exhaustive, is also useful to illustrate some of the characteristics of our methods. We remark that the properties of Algorithms \ref{algoCompl} and \ref{algoComplb} have been studied in the previous section under the assumption that $K$ is bounded. But it is easy to see from the proofs that this condition can be substituted by any of the assumptions studied in this section that still guarantee boundedness of the sequence generated by Algorithms \ref{algoCompl}
and \ref{algoComplb}. The only adjustment that needs to be made is that the Lipschitz constants  used in the proofs of Theorems \ref{th:complbad} and \ref{th:complbadb} are no longer the Lipschitz constants on $K$ but rather the ones on the compact set $S$, which  is shown to contain the sequence $\{x^\nu\}$, and that any reference to the boundedness of $K$ should be substituted by a reference to the compactness of $S$. \textcolor{black}{Of course, in order for this approach to be sensible as a complexity bound, the set $S$ must be determined a priori and should not depend on the sequence generated by the algorithm.}

 1. [Valid for Algorithms \ref{algoBasic}, \ref{algoCompl} and \ref{algoComplb}] The boundedness of $K$
 obviously guarantees the boundedness of $\{x^\nu\}$ for all the algorithms we considered. We already used this fact for Algorithms \ref{algoCompl} and \ref{algoComplb}, but the same result holds also for
 Algorithm \ref{algoBasic}: we report this case here for completeness and uniformity of presentation.
 In fact,
 for all the three algorithms we have that $\gamma^\nu \in (0,1]$ and therefore the constraint $d\in K-x^\nu$ in \eqref{eq:p_k} and
 the convexity of $K$ guarantee that if $x^\nu\in K$ then also
 $x^\nu + \gamma^\nu d(x^\nu)$ belongs to $K$.
 This case covers most instances of practical interest since, in basically all real-world problems, variables are
 naturally limited by lower and upper bounds and we can take $K$ to be  the rectangle defined by these
 quantities.

 2. [Valid for Algorithms \ref{algoBasic}, \ref{algoCompl} and \ref{algoComplb}] Another setting in which we can guarantee the boundedness of the iterations is when we turn our schemes into  feasible  methods by choosing $\tilde g_i$s that are UCAs of the $g_i$s and a feasible starting point $x^0$, see Section \ref{sec:subsec3.1}.

2a.
 In this setting, assume that the following classical condition holds:
\begin{equation}\label{eq:ass bound1}
{\cal L_1} \, \triangleq \,  \{x\in K : g(x) \,\leq\, 0,\, f(x)\, \leq\, f(x^0)\} \;\, \;\mbox{\rm is bounded},
\end{equation}
i.e., the level set of value $f(x^0)$ for the objective function intersected with the feasible set is bounded. Then,
if we also assume that $\tilde f$ is an UCA of $f$, see Remark \ref{rem:feasible DSM}, we can show that
the sequence $\{x^\nu\}$ generated by any  of the Algorithms \ref{algoBasic}, \ref{algoCompl} and
\ref{algoComplb} is contained in the bounded
set $\cal L_1$. Since the sequence $\{x^\nu\}$ belongs to $\mathcal X$,  it is enough to show that, at each
iteration,
$f(x^{\nu + 1}) \le f(x^\nu)$. To this end, observe that, since each $x^\nu$ is feasible, we always have $\max_i
\{g_i(x^\nu)_+\}=\kappa(x^\nu)=\theta(x^\nu)=0$, and $d =0$ is feasible for (P$_{x^\nu}$).
Therefore, applying the minimum principle to (P$_{x^\nu}$), we have $\nabla_1 \tilde f(d(x^\nu); x^\nu)\trt (0 -
d(x^\nu)) \, \geq\, 0$ and, in turn
\begin{equation}\label{eq:minimum principle}
\nabla_1 \tilde f(d(x^\nu); x^\nu)\trt d(x^\nu) \, \leq\, 0.
\end{equation}
Since $\tilde f$ is (strongly) convex, we get
\begin{equation}\label{eq:deriv1}
\begin{array}{rcl}
f(x^\nu) &= & \tilde f(0;x^\nu) \, \geq \, \tilde f(\gamma^\nu d(x^\nu); x^\nu) +
\nabla _1 \tilde f (\gamma^\nu d(x^\nu); x^\nu)\trt(0 -\gamma^\nu d(x^\nu))
\\[0.5em]
& \geq & f(x^\nu + \gamma^\nu d(x^\nu)) -\gamma^\nu \nabla \tilde f (\gamma^\nu d(x^\nu); x^\nu) \trt
d(x^\nu),
\end{array}
\end{equation}
where the second inequality follows from \eqref{eq:uppera}.
By the strong convexity with modulus $c$ of $\tilde f$, see A1, we can also write
\[
\left(\nabla_1 \tilde f(d(x^\nu);x^\nu) - \nabla_1 \tilde f(\gamma^\nu d(x^\nu);x^\nu)\right)\trt (d(x^\nu) -
\gamma^\nu d(x^\nu)) \, \geq\,
c (1-\gamma^\nu)^2 \|d(x^\nu)\|^2,
\]
which, with simple manipulations, yields
\[
-\gamma^\nu \nabla_1 \tilde f(\gamma^\nu d(x^\nu);x^\nu)\trt d(x^\nu) \, \geq \,
-\gamma^\nu \nabla_1 \tilde f(d(x^\nu);x^\nu)\trt d(x^\nu) + c\gamma^\nu (1-\gamma^\nu) \|d(x^\nu)\|^2.
\]
 Plugging this inequality in \eqref{eq:deriv1} we get
 \[
 f(x^\nu + \gamma^\nu d(x^\nu)) \, \leq \,f(x^\nu) + \gamma^\nu \nabla_1 \tilde f(d(x^\nu); x^\nu)\trt d(x^
 \nu)
 \]
 which, in view of \eqref{eq:minimum principle}, shows that $f(x^{\nu+1})\leq f(x^\nu)$ as desired for any
 choice of $\gamma^\nu$.

2b.  The requirement that $\tilde f$  be an upper approximation of $f$ is actually not needed for Algorithm \ref{algoComplb}. In fact, noting that  $W(x^\nu;\varepsilon) = f(x^\nu)$ for any feasible $x^\nu$ and for any positive $\varepsilon$, step \ref{S.7b} in Algorithm \ref{algoComplb} guarantees $f(x^{\nu+1}) \leq f(x^\nu)$. In fact, in the current setting, a suitable stepsize can be found in step \ref{S.7b} because of \eqref{eq:condbound} (that still holds even if $K$ is not bounded, see the conditions in the if-block at step \ref{S.3b}), recalling that $\theta(x^\nu)=0$ for any feasible $x^\nu$. It is therefore clear that the whole sequence $\{x^\nu\}$ is contained in  ${\cal L}_1$.

2c. We can avoid the UCA requirement on  $\tilde f$ also for Algorithm \ref{algoCompl}, provided we assume that  $\nabla f$ is Lipschitz continuous on  $K$  (with modulus $L_{\nabla f}$). By the descent lemma we can write
\[
f(x^\nu + \gamma^\nu d(x^\nu)) - f(x^\nu) \leq \gamma^\nu \nabla f(x^\nu)\trt d(x^\nu)
+ (\gamma^\nu)^2\frac{L_{\nabla f}}{2}\|d(x^\nu)\|^2.
\]
In turn, we get from \eqref{eq:preldisc2} (which, again, is easily seen to hold in the current setting), taking into account that $\theta(x^\nu) =0$ because $x^\nu$ is feasible,
\begin{equation}\label{feasdescf}
\displaystyle
f(x^\nu + \gamma^\nu d(x^\nu)) - f(x^\nu) \le -\gamma^\nu \left(\eta c - \frac{\gamma^\nu}{2}L_{\nabla f} \right)
\|d(x^\nu)\|^2,
\end{equation}
for every $\nu$. The instructions in {\bf Data} of Algorithm \ref{algoCompl} are easily seen to entail $\gamma^\nu \in (0, \min\{1, {2 \eta c}/{L_{\nabla f}}\}]$, so that \eqref{feasdescf} implies that $\{x^\nu\}$ is all contained in ${\cal L}_1$, and therefore that $\{x^\nu\}$ is bounded if ${\cal L}_1$ is bounded.

2d. If we want to eliminate the UCA property of $\tilde f$ also for Algorithm \ref{algoBasic}, we need again $\nabla f$ to be Lipschitz continuous on $K$ and to strengthen condition \eqref{eq:ass bound1}, requiring that
 \begin{equation}\label{eq:ass bound1nof}
 {\cal L_2^\alpha} \triangleq \{x \in K : g(x) \le 0, f(x) \le \alpha\}\;\, \;\mbox{\rm is bounded for every $\alpha
 \in\mathbb R$}.
 \end{equation}
 Then, invoking \eqref{eq:convscfinter}, as soon as $\gamma^\nu$ becomes smaller than $\frac{2c}{L_{\nabla f}}$, we stay in the set  $ {\cal L_2^\alpha} $ for some value of $\alpha$; by \eqref{eq:ass bound1nof} this implies the boundedness of $\{x^\nu\}$.

2e. Finally, it is worth observing that if the feasible set $\cal X$ is bounded, the use of UCAs for the constraints $g$ is enough to guarantee the boundedness of $\{x^k\}$ since, if we start with a feasible point, $\{x^\nu\}$ remains feasible whatever the algorithm we use, see Section \ref{sec:subsec3.1}.

3. [Valid for Algorithms \ref{algoBasic}, \ref{algoCompl} and \ref{algoComplb}]
 Knowing a feasible point $x^0$ to start the algorithm from can be difficult in some applications. But fortunately,  the results in point 2 above can be generalized in order to avoid the feasibility requirement.

3a. Suppose that we start the algorithm with a possibly infeasible point $x^0 \in K$. Assume that
\begin{equation}\label{eq:ass bound2}
{\cal L}_3\, \triangleq\,  \{x\in K : g(x) \,\leq\, \max_i\{ g_i(x^0)_+\}\}\;\, \;\mbox{\rm is bounded}.
\end{equation}
If we use $\tilde g_i$s that are UCAs for the $g_i$s, we can show by induction that the whole sequence $\{x^\nu\}$ generated by  Algorithms \ref{algoBasic}, \ref{algoCompl} and \ref{algoComplb} is contained in ${\cal L}_3 $. In fact, the starting point $x^0$ of course belongs to ${\cal L}_3 $. Suppose now that  $x^\nu \in {\cal L}_3 $, meaning that
 \[\tilde g_i(0;x^\nu) = g_i(x^\nu)  \leq \max_i\{ g_i(x^0)_+\}.\]
 We also have
 \[
\tilde g_i(d(x^\nu); x^\nu) \,\leq \,\kappa(x^\nu)\, \leq\,
  \max_i\{ g_i(x^\nu)_+\},
 \]
 where the first inequality is just feasibility for subproblem (P$_{x^\nu}$), and the second one follows by the definition of $\kappa(x^\nu)$.  The last two displayed formulas show that, for any $\gamma^\nu \in [0,1]$,
\[
g_i(x^\nu + \gamma^\nu d(x^\nu)) \leq \tilde g_i  (\gamma^\nu d(x^\nu); x^\nu)
\leq \max_i\{ g_i(x^\nu)_+\},
 \]
 where the first inequality is due to the UCA property \eqref{eq:upperb}, while the second relation derives from the convexity of $g_i(\cdot;x^\nu)$.

  3b.
When the constraints $g_i$s are convex, it is well known that the boundedness of ${\cal L}_3 $ holds if and only if the feasible set is bounded. Then, in principle we can set $\tilde g (d;x) = g(x+d)$ and only approximate the objective function. This particular $\tilde g$ is a UCA, indeed. This approach seems particularly well suited to the case in which the $g_i$s are linear because the resulting subproblem then has simple linear constraints. Note also that keeping the original (convex) constraints in the subproblems is something routinely done in most MM methods.

4. [Valid for Algorithms \ref{algoBasic}, \ref{algoCompl} and \ref{algoComplb}]
Another interesting case arises if we suppose that the eMFCQ holds  and
\begin{equation}\label{eq:ass bound3}
{\cal L}_4^\alpha \, \triangleq\,
\{x\in K : \max_i\{ g_i(x)_+\} \le \alpha\} \;\,\; \mbox{\rm is bounded for every $\alpha \in \mathbb R_+$.}
\end{equation} Note that \eqref{eq:ass bound3} simply states that
the function $\max_i\{ g_i(x)_+\}$ is coercive. We can
therefore  find
positive $\alpha_1$ and $\alpha_2$ such that if
 $x^\nu\in {\cal L}_4^{\alpha_1}$ then $x^\nu + \gamma^\nu d(x^\nu) \in {\cal L}_4^{\alpha_2}$ for all
 $\gamma^\nu \in (0,1]$.
 Now,  following the same line of reasoning as relations \eqref{eq:nonsmmaj1ter}, we have
\begin{equation}\label{feasdesc}
\displaystyle
\max_i\{ g_i(x^\nu + \gamma^\nu d(x^\nu))_+\} - \max_i\{ g_i(x^\nu)_+\}  \le  -\gamma^\nu \left(
\theta(x^\nu) - \frac{\gamma^\nu}{2} \max_i \{L_{\nabla g_i}\} \|d(x^\nu)\|^2\right)
\end{equation}
for every $x^\nu\in  {\cal L}_4^{\alpha_2}$, where $L_{\nabla g_i}$ are Lipschitz constants of the gradients of
$g_i$ on  ${\cal L}_4^{\alpha_2}$; we remark that since  ${\cal L}_4^{\alpha_2}$ is bounded, existence of these
constants is a very mild requirement. Denote by $\bar \theta >0$ a positive constant such that
$\theta(x) \geq \bar \theta$ for all points in the set $\Delta \triangleq {\cal L}_4^{\alpha_2}\setminus \int {\cal
L}_4^{\alpha_1}$; note that this set is compact by \eqref{eq:ass bound3}. Such $\bar \theta$ surely exists
because  the eMFCQ implies there are no ES in the  set $\Delta$ and therefore the continuous function $\theta(x)
$ is positive on $\Delta$. By \eqref{feasdesc}, we can then write, for any $x^\nu \in \Delta$,
\begin{equation}\label{feasdesc1}
\displaystyle
\max_i\{ g_i(x^\nu + \gamma^\nu d(x^\nu))_+\} - \max_i\{ g_i(x^\nu)_+\}  \le  -\gamma^\nu \left(\bar \theta
- \frac{\gamma^\nu}{2} \max_i \{L_{\nabla g_i}\} \beta^2\right).
\end{equation}
It is then clear that a threshold value $\bar \gamma>0$ exists such that, if $\gamma^\nu \leq \bar \gamma$,
then $\max_i\{ g_i(x^\nu + \gamma^\nu d(x^\nu))_+\} \leq \max_i\{ g_i(x^\nu)_+\}  $. Now, two cases can
occur. If $x^\nu$ belongs to $\int {\cal L}_4^{\alpha_1}$, then, by how we have chosen $\alpha_2$, $x^{\nu+1}
$ belongs to ${\cal L}_4^{\alpha_2}$. If instead  $x^\nu$  belongs to  $\Delta$, by taking $\gamma^\nu \leq
\bar \gamma$, we are again sure that $x^{\nu+1}$ still belongs to ${\cal L}_4^{\alpha_2}$. We can so conclude
that by using stepsizes smaller that $\bar \gamma$, iterations never leave the set ${\cal L}_4^{\alpha_2}$ and
therefore stay bounded.

5. [Valid for Algorithms \ref{algoCompl} and \ref{algoComplb}]
The eMFCQ assumption in case 4 can be replaced by the requirement that $f$ be bounded from below on $K$ if one employs Algorithm \ref{algoComplb}, a condition to which the Lipschitz continuity of $\nabla f$ and $\nabla g_i$ on $K$ has to be added when Algorithm \ref{algoCompl} is resorted to. For both cases, the proof of the claim reduces to showing that, even without requiring $K$ to be bounded, we still have the sufficient descent condition
\begin{equation}\label{eq:asdesc23}
	W(x^{\nu+1}; T^\nu) - W(x^\nu; T^\nu) \le -\gamma^\nu \frac{\eta c}{4} \|d(x^\nu)\|^2
\end{equation}
for every $\nu$. In fact, once relation \eqref{eq:asdesc23} has been proven to be valid, in turn we get \[
\begin{array}{rcl}
T^{\nu+1} (f(x^{\nu+1})-\bar f) + \max_i \{g_i(x^{\nu+1})_+\} & \le & T^{\nu} (f(x^{\nu+1})-\bar f) + \max_i
\{g_i(x^{\nu+1})_+\}\\[5pt]
& \le & T^\nu (f(x^\nu) - \bar f) + \max_i \{g_i(x^\nu)_+\} - T^\nu \gamma^\nu \frac{\eta c}{4} \|d(x^\nu)\|
^2,
\end{array}
\]
where the first relation follows from observing that $T^\nu$ is non increasing. The sequence generated by the algorithms is now easily shown to be bounded. Indeed, the inequality above shows that the nonnegative sequence $\{T^\nu (f(x^\nu) - \bar f) + \max_i \{g_i(x^\nu)_+\}\}$ is non increasing and therefore convergent. Suppose now by contradiction that $\{x^k\}$ is unbounded. By \eqref{eq:ass bound3} this implies that $\max_i \{g_i(x^\nu)_+\}$ goes to infinity; since $f(x^\nu) - \bar f$ is nonnegative, this contradicts the convergence of the sequence $\{T^\nu (f(x^\nu) - \bar f) + \max_i \{g_i(x^\nu)_+\}\}$.

Let us now show why, in the current setting, \eqref{eq:asdesc23} is still satisfied for Algorithms \ref{algoCompl} and \ref{algoComplb}.

5a.
Concerning Algorithm \ref{algoComplb}, \eqref{eq:asdesc23}
is enforced as the algorithm progresses, see \ref{S.7b}. We remark that, even in the present setting, given an iterate $x^\nu$, in \ref{S.7b} a sufficiently small stepsize $\gamma^\nu$ still exists such that \eqref{eq:asdesc23} is verified: this follows by standard reasoning ab absurdo in view of \eqref{eq:condbound}, which still holds even if $K$ is not bounded (see the conditions in the if-block at step \ref{S.3b}), and observing that the directional derivative of $\max_i \{g_i(x^\nu)_+\}$ is bounded from above by $-\theta(x^\nu)$, thanks to \eqref{eq:convscgter}.

5b.
As for Algorithm \ref{algoCompl}, with $\nabla f$ and $\nabla g_i$ assumed to be Lipschitz continuous on $K$, let $L_{\nabla f}$ and $L_{\nabla g_i}$ be the corresponding Lipschitz moduli. Again, without requiring $K$ to be bounded as done in Section \ref{Sec:complexity}, condition \eqref{eq:asdesc23} is clearly satisfied, since \eqref{eq:preldisc} and \eqref{eq:preldisc2} remain valid following the same line of reasoning as relation \eqref{eq:nonsmmaj1ter} (here with $K$ not assumed to be bounded, but under the Lipschitz continuity of $\nabla f$ and $\nabla g$), and as a straightforward consequence of the conditions in the if-block at step \ref{S.3}, respectively.

Note that in both cases 5a and 5b, although $\{x^\nu\}$ is contained in ${\cal L}_4^\alpha$ for some $\alpha$, this quantity is possibly unknown in advance.

\smallskip
We summarize the above  conditions implying boundedness in Table \ref{tab2}. We also clarify (see the last column that only applies to Algorithms \ref{algoCompl} and \ref{algoComplb}) on a case by case basis if these assumptions make it possible to perform an Iteration Complexity (IC) or a Global Convergence Rate (GCR) analysis.
In fact, when the bounded set to which the sequence $\{x^\nu\}$ belongs is defined by means of quantities that are known in advance, for example if it is ${\cal L}_1$, we can still speak of iteration complexity results derived for Algorithms \ref{algoCompl} and \ref{algoComplb}, as done in Theorems \ref{th:complbad} and \ref{th:complbadb}; when the set is known to exist, but is not known beforehand (as for example in cases 5a and 5b), we obtain instead a global convergence rate for the corresponding algorithms, since the constants involved in the big ${\mathcal O}$ bounds cannot be determined a priori.


\begin{table}[h]
\TABLE
{Summary of conditions for boundedness of iterates \label{tab2}}
{\begin{tabular}{c|c|c|c|c|c|c|c}
& $K$ & $\tilde f$ &  $\tilde g$ &  $x^0$ &  Other assumptions & Algorithm & IC or GCR\\ \hline
1 & bounded & - & - &-& -& \ref{algoBasic}, \ref{algoCompl}, \ref{algoComplb} & IC \\[5pt]
2a & -& UCA \eqref{eq:uppera}& UCA \eqref{eq:upperb}&feasible &  ${\cal L_1} $ bounded, see \eqref{eq:ass
bound1} & \ref{algoBasic}, \ref{algoCompl}, \ref{algoComplb} & IC \\[5pt]
2b & -& - & UCA \eqref{eq:upperb}&feasible &  ${\cal L_1} $ bounded, see \eqref{eq:ass bound1} &
\ref{algoComplb} & IC\\[5pt]
2c & -& - & UCA \eqref{eq:upperb}&feasible &  ${\cal L_1} $ bounded, see \eqref{eq:ass bound1}, $\nabla f$
Lipschitz & \ref{algoCompl}, \ref{algoComplb} & IC\\[5pt]
2d & -& - & UCA \eqref{eq:upperb}&feasible &  ${\cal L_2^\alpha} $ bounded, see \eqref{eq:ass bound1nof}, $
\nabla f$ Lipschitz & \ref{algoBasic}, \ref{algoCompl}, \ref{algoComplb}& IC\\[5pt]
2e & -& - & UCA \eqref{eq:upperb}&feasible & $\mathcal X$ bounded & \ref{algoBasic}, \ref{algoCompl},
\ref{algoComplb} & IC\\[5pt]
3a & -& - & UCA \eqref{eq:upperb}& -  &  ${\cal L}_3$ bounded, see \eqref{eq:ass bound2} & \ref{algoBasic},
\ref{algoCompl}, \ref{algoComplb} & IC \\[0.5em]
3b & -& - & $\tilde g = g$ & - &  $g_i$s convex, $\mathcal X$ bounded & \ref{algoBasic}, \ref{algoCompl},
\ref{algoComplb} & IC\\[5pt]
4 & -& - & -  & -  &  ${\cal L}_4^\alpha $ bounded, see \eqref{eq:ass bound3}, eMFCQ & \ref{algoBasic},
\ref{algoCompl}, \ref{algoComplb} & GCR\\[5pt]
5a &-&-&-&-& ${\cal L}_4^\alpha $ bounded, see \eqref{eq:ass bound3}, $f$ low. bounded & \ref{algoComplb} & GCR\\[5pt]
5b &-&-&-&-& ${\cal L}_4^\alpha $ bounded, see \eqref{eq:ass bound3}, $f$ low. bounded,& \ref{algoCompl}, \ref{algoComplb}  & GCR\\
 &&&&&$\nabla f$ and $\nabla g$ Lipschitz &
 \end{tabular}}{}
 \end{table}

\section*{Acknowledgments.} The authors are very thankful to Philippe Toint, whose detailed comments on an earlier version of this paper helped shape it. The constructive criticisms of two Referees were also very useful in improving the quality of the paper.

Facchinei was partially supported by MIUR  PLATINO (PLATform for INnOvative services in future internet) PON project, under Grant Agreement no. PON01\_01007.
Scutari was supported by the USA National Science Foundation (NSF) under Grants CIF
1564044    and CAREER Award No. 1555850, and the Office of Naval Research (ONR)
Grant N00014-16-1-2244. Kungurtsev was supported by the Czech Science Foundation project 17-26999S and the OP VVV project
CZ.02.1.01/0.0/0.0/16\_019/0000765 Research Center for Informatics


\bibliographystyle{informs2014}
\bibliography{Surbib}



\end{document}